\documentclass[11pt]{amsart}
\usepackage{amsfonts}

\newcommand{\Z}{{\mathbb{Z}}}

\newcommand{\C}{{\mathbb{C}}}

\newcommand\modA{\operatorname{-mod}}

\newcommand{\GL}{{\operatorname{GL}}}
\newcommand{\Aut}{{\operatorname{Aut}}}

\newcommand{\ulambda}{{\boldsymbol{\lambda}}}
\newcommand{\umu}{{\boldsymbol{\mu}}}

\newcommand{\Fock}{\mathcal F}
\newcommand{\residue}{\mathrm res}
\newtheorem{thm}{Theorem}[section]
\newtheorem{lem}[thm]{Lemma}
\newtheorem{cor}[thm]{Corollary}
\newtheorem{prop}[thm]{Proposition}

\theoremstyle{definition}
\newtheorem{exmp}[thm]{Example}
\newtheorem{defn}[thm]{Definition}

\theoremstyle{remark}

\newcommand{\Fl}{\mathcal{F}\ell}
\newcommand{\Nil}{\mathcal{N}}

\newcommand{\Id}{\operatorname{Id}}
\newcommand{\Ker}{\operatorname{Ker}}
\newcommand{\Rad}{\operatorname{Rad}}
\newcommand{\Top}{\operatorname{Top}}
\newcommand{\Soc}{\operatorname{Soc}}
\newcommand{\Ind}{\operatorname{Ind}}
\newcommand{\Res}{\operatorname{Res}}
\newcommand{\Ext}{\operatorname{Ext}}
\newcommand{\Hom}{\operatorname{Hom}}
\newcommand{\Endhom}{\operatorname{End}}
\newcommand{\wt}{\operatorname{wt}}
\newcommand{\diag}{\operatorname{diag}}
\newcommand{\supp}{\operatorname{supp}}
\newcommand{\Sym}{\mathfrak S}
\newcommand{\Hecke}{H}

\begin{document}
\title{The modular branching rule for affine Hecke algebras of type $A$}
\author{Susumu Ariki, Nicolas Jacon and C\'edric Lecouvey}
\address{S.A.: Research Institute for Mathematical Sciences, Kyoto University, Kyoto
606-8502, Japan}
\email{ariki@kurims.kyoto-u.ac.jp}
\address{N.J.:Universit\'e de Franche-Comt\'e, UFR Sciences et Techniques, 16 route
de Gray, 25 030 Besan\c{c}on, France.}
\email{njacon@univ-fcomte.fr}
\address{C. L.: Laboratoire de Math\'{e}matiques Pures et Appliqu\'{e}es Joseph
Liouville Centre Universitaire de la Mi-Voix B.P. 699 62228 Calais France}
\email{Cedric.Lecouvey@lmpa.univ-littoral.fr}
\date{Revised Mar.19, 2010}

\begin{abstract}
For the affine Hecke algebra of type $A$ at roots of unity, 
we make explicit the correspondence between 
geometrically constructed simple modules and combinatorially 
constructed simple modules and prove the modular branching rule. 
The latter generalizes work by Vazirani.
\end{abstract}

\maketitle


\pagestyle{myheadings}

\markboth{Susumu Ariki, Nicolas Jacon and C\'edric Lecouvey}{Affine Hecke
algebras of type $A$}


\section{Introduction}

In \cite{CG}, Ginzburg explains his geometric construction of simple modules
over (extended) affine Hecke algebras $\Hecke_n$ defined over $\C$. In this paper, 
we consider the affine Hecke algebra of type $A$ whose parameter is a root of unity. 
Then, the simple modules are labelled by aperiodic multisegments. 

On the other hand, Dipper, James and Mathas' Specht module theory gives us a
combinatorial construction of simple modules of cyclotomic Hecke algebras, 
and they exhaust all the simple modules of the affine Hecke algebra. 
The simple modules are labelled by Kleshchev multipartitions. 

If one wants to compute something about simples, 
the combinatorially defined simple modules often have more advantage than 
the geometrically defined simple modules, and we may work over 
any algebraically closed field. On the other hand, 
the geometrically defined simple modules are very useful 
in several circumstances. Hence, explicit description of 
the module correspondence between the two constructions is desirable. 

We provide this explicit description of the module 
correspondence in this article. Note that 
both the set of aperiodic multisegments and the set of Kleshchev multipartitions 
have structure of Kashiwara crystals. Then, we show that 
the crystal embedding gives the module correspondence. We also describe the 
crystal embedding explicitly. 

Closely related to this result is the modular branching rule. One may prove 
the result on the module correspondence by using this, which is our first proof, 
or one may prove the modular branching rule by first establishing the result 
on the module correspondence, which is our second proof. Note that we mean here 
the modular branching rule in the original sense as we explain in the next paragraph. 
Some authors use the terminology in weaker sense, which does not imply the module correspondence. 

Let $L_\psi$ be the simple module labelled by a multisegment $\psi$, 
whose precise meaning will be explained in section 4. 
The \emph{modular branching rule} is a rule to describe 
$\Soc(i\text{-}\Res_{\Hecke_{n-1}}^{\Hecke_n}(L_\psi))$, or
equivalently $\Top(i\text{-}\Res_{\Hecke_{n-1}}^{\Hecke_n}(L_\psi))$. We show that 
$$
\Soc(i\text{-}\Res_{\Hecke_{n-1}}^{\Hecke_n}(L_\psi))
=L_{\tilde e_i\psi},
$$
where $\tilde e_i$ is the Kashiwara operator. We give a geometric proof of this 
rule in the framework of Lusztig and Ginzburg's theory. This gives the first proof. 
On the other hand, if one uses results in \cite{Ari3} and \cite{Ari4}, 
both become easier, and this is the second proof. 

Recall that the main result of \cite{V} is the modular
branching rule when the parameter of the affine Hecke algebra is not a root of unity. 
Hence our result generalizes \cite[Theorem 3.1]{V}. In \cite{Ari2} and \cite{AM} it was proved that 
affine $sl_e$ controls the modular representation theory of cyclotomic 
Hecke algebras. Later\footnote{Compare the submission date of \cite{AM} with \cite{G}.}, 
Grojnowski gave another proof \cite[Theorem 14.2, 14.3]{G} for several main results in \cite{AM}. 
We note here that he writes in \cite[14.1]{G} that 
his IMRN paper proved the part concerning the canonical basis in \cite{Ari2}, but 
his announcement in 1995 was that he had some computation of Kazhdan-Lusztig polynomials 
which he was able to deduce from the IMRN paper: 
the use of $i$-induction and $i$-restriction functors, integrable modules over 
$\mathfrak g(A^{(1)}_{e-1})$, and Lusztig's aperiodicity were absent in the assertion. 

The idea of his proof in \cite{G} came from Leclerc's observation that Kleshchev and Brundan's work on 
the modular branching rule of the symmetric group and the Hecke algebra of type $A$ may be 
understood in crystal language. The proof is interesting, but the modular branching 
rule in our sense is not proved in \cite{G} and it is natural to ask whether the crystal 
he used coincides with the one used in \cite{Ari2} and \cite{AM}. 
It was settled affirmatively in \cite{Ari4}, but it still used several results from \cite{G}. 
Here in this paper, the modular branching rule for the affine Hecke algebra, which is a stronger 
statement than the statement in \cite{G} that the modular branching gives a crystal which is 
isomorphic to $B(\infty)$, is proved in a direct manner. 
It still uses the multiplicity one result from \cite{GV}, but it replaces \cite{G}.

The paper is organized as follows. In section 2, we review basic facts on 
the crystal $B(\infty)$ of type $A^{(1)}_{e-1}$. In section 3, we prepare 
for a geometric proof of the modular branching 
rule of the affine Hecke algebra. 
In section 4, we give the geometric proof of the modular branching rule. 
In section 5, we introduce crystals of deformed Fock spaces and state results 
to compute crystal isomorphisms among them. 
In section 6, we prove a lemma on the module correspondence of simple modules in 
various labellings and 
give a combinatorial proof of the modular branching rule in the framework 
of Fock space theory for cyclotomic Hecke algebras. \\
\\
\textbf{Acknowledgements.}  Part of this work was done while 
the authors were visiting the MSRI in Berkeley in 2008. 
The authors wish to thank the institute for its hospitality 
and the organizers of the two programs for their invitation. 
The second author is also grateful to Hyohe Miyachi for 
fruitful discussions there. 
The second author is supported by the ``Agence Nationale de la~Recherche"  
(project JCO7-192339).

\section{Preliminaries}

Let $e\ge2$ be a fixed integer, $\mathfrak g$ the Kac-Moody Lie algebra of 
type $A^{(1)}_{e-1}$. We denote by $U_v^-$ the negative part 
of the quantum affine algebra $U_v(\mathfrak g)$, which is generated by the 
Chevalley generators $f_i$, where $i\in\Z/e\Z$, subject to the quantum 
Serre relations. In this section, we review basic facts on $U_v^-$ 
and its crystal. We denote the simple roots by $\alpha_i$, 
and the simple coroots by $\alpha_i^\vee$, for $i\in\Z/e\Z$.

\subsection{The crystal $B(\infty)$}
Let us introduce the Kashiwara operator $\tilde f_i$, for $i\in\Z/e\Z$, 
on $U_v^-$. Let $e_i$, $i\in\Z/e\Z$, 
be Chevalley generators of the positive part of $U_v(\mathfrak g)$
and $t_i=v^{\alpha_i^\vee}$. The following two lemmas are due to Kashiwara. 

\begin{lem}
For each $u\in U_v^-$, there exist unique $u'$ and $u''$ in $U_v^-$ such that 
we have
$$
e_iu-ue_i=\frac{t_iu'-t_i^{-1}u''}{v-v^{-1}}.
$$
\end{lem}

We define an operator $e_i'$ on $U_v^-$ by $e_i'u=u''$, for $u\in U_v^-$. 
The algebra generated by $\{f_i\}_{i\in\Z/e\Z}$ and 
$\{e_i'\}_{i\in\Z/e\Z}$ is called the \emph{Kashiwara algebra}. 
Let $f_i^{(n)}$ be the $n^{th}$ divided power of $f_i$. 

\begin{lem}
Let $P\in U_v^-$. For each $i\in\Z/e\Z$, we may write 
$$
P=\sum_{n\in\Z_{\ge0}}f_i^{(n)}u_n, 
$$
where $u_n\in U_v^-$ are such that $e_i'u_n=0$, for all $n$. Further, 
the expression is unique. 
\end{lem}

We define $\tilde e_iP=\sum_{n\in\Z_{\ge1}}f_i^{(n-1)}u_n$
and $\tilde f_iP=\sum_{n\in\Z_{\ge0}}f_i^{(n+1)}u_n$. They are well-defined. 
Let $R$ be the subring of $\C(v)$ consisting of elements which are regular at 
$v=0$. Then, we define 
$$
L(\infty)=\sum_{N\in\Z_{\ge0}}
\sum_{(i_1,\dots,i_N)\in(\Z/e\Z)^N} R\tilde f_{i_1}\cdots\tilde f_{i_N}1
$$
and 
$$
B(\infty)=\left(\cup_{N\in\Z_{\ge0}}\cup_{(i_1,\dots,i_N)\in(\Z/e\Z)^N}
\tilde f_{i_1}\cdots\tilde f_{i_N}1+vL(\infty)\right)\setminus\{0\}.
$$
$B(\infty)$ is a basis of the $\C$-vector space $L(\infty)/vL(\infty)$. 
$U_v^-$ admits a root space decomposition
$U_v^-=\oplus_{\alpha\in Q_+} (U_v^-)_{-\alpha}$, where 
$Q_+=\sum_{i\in\Z/e\Z}\Z_{\ge0}\alpha_i$, and it follows that 
$$
B(\infty)=\bigsqcup_{\alpha\in Q_+} B(\infty)_{-\alpha}.
$$
We define $\wt(b)=-\alpha$ if $b\in B(\infty)_{-\alpha}$. 
Then, by defining 
$$
\epsilon_i(b)=\max\{k\in\Z_{\ge0} \mid \tilde e_i^kb\ne0\}\;\;\text{and}\;\; 
\varphi_i(b)=\epsilon_i(b)+\wt(b)(\alpha_i^\vee),
$$
for $b\in B(\infty)$, 
$(B(\infty), \wt, \epsilon_i, \varphi_i, \tilde e_i, \tilde f_i)$ is a 
$\mathfrak g$-crystal in the sense of Kashiwara \cite[p.48]{Ka}. 

We define the bar operation on $U_v^-$ by $\bar v=v^{-1}$ and $\bar f_i=f_i$. 
Lusztig and Kashiwara independently constructed the canonical basis/the global basis 
$$
\{G_v(b) \mid b\in B(\infty)\}
$$
of $U_v^-$, which is characterized by the property that 
$$
\overline{G_v(b)}=G_v(b),\quad G_v(b)+vL(\infty)=b.
$$

\begin{exmp}
Let $e=3$. Then, $e_2$ and $f_1$ commute so that $e_2'f_1=0$ and
$\tilde f_2f_1=f_2f_1$ follows. Similarly, $\tilde f_1f_2=f_1f_2$. 
Thus, $\{f_1f_2, f_2f_1\}$ is the canonical basis of 
$(U_v^-)_{-\alpha_1-\alpha_2}$. For the null root 
$\delta=\alpha_0+\alpha_1+\alpha_2$, 
$\{f_0f_1f_2, f_0f_2f_1, f_1f_0f_2, f_1f_2f_0, f_2f_0f_1, f_2f_1f_0\}$ is the canonical 
basis of $(U_v^-)_{-\delta}$. Of course, more complex linear combinations of 
monomials in $f_i$ appear in the canonical basis 
of other $(U_v^-)_{-\alpha}$. 
\end{exmp}

\subsection{Hall algebras}

The crystal $B(\infty)$ has a concrete description. 
Let $\Gamma_e$ be the cyclic quiver of length $e$. 
This is an oriented graph with vertices $\Z/e\Z$ and 
edges $\{(i,i+1),i\in\Z/e\Z\}$. Let $V=\oplus_{i\in\Z/e\Z} V_i$ 
be a finite dimensional $\Z/e\Z$-graded vector space, and define
$$
E_V=\bigoplus_{i\in\Z/e\Z}\Hom_\C(V_i, V_{i+1})\subseteq \Endhom_\C(V).
$$
An element $X\in E_V$ is called a \emph{representation} of $\Gamma_e$ on $V$. 
The vector 
$$
\underline\dim\;V=(\dim V_i)_{i\in\Z/e\Z}
$$
is called the \emph{dimension vector} of the representation. 

If $V$ runs through all finite dimensional $\Z/e\Z$-graded vector spaces, we 
obtain the category of representations of $\Gamma_e$. 
It is the same as the category of finite dimensional $\C\Gamma_e$-modules, 
where $\C\Gamma_e$ is the path algebra of $\Gamma_e$. 
If $X$ is nilpotent as an endomorphism of $V$, we say that 
the representation $X$, or the corresponding $\C\Gamma_e$-module, is 
\emph{nilpotent}. We denote by $\Nil_V$ 
the subset of nilpotent representations in $E_V$. 
Let $G_V=\prod_{i\in\Z/e\Z}\GL(V_i)$. It acts on $E_V$ and $\Nil_V$ by 
conjugation and two representations are equivalent if and only if they are in the same $G_V$-orbit. 

For each $i\in \Z/e\Z$, let $V=V_i=\C$ and $X=0$. Then it defines a simple 
$\C\Gamma_e$-module. 
We denote it by $S_i$. They are nilpotent representations. 

\begin{exmp}
Let $G_n=\GL_n(\C)$ and 
suppose that $s\in G_n$ has order $e$. Let $\zeta$ 
be a primitive $e^{th}$ root of unity, $V=\C^n$, 
and let $V_i$ be the eigenspace of $s$ for the eigenvalue $\zeta^i$. 
If $X\in \Endhom_\C(V)$ is such that $sXs^{-1}=\zeta X$ then  
$XV_i\subseteq V_{i+1}$. 
Thus, $X$ defines a representation of $\Gamma_e$ on $V$. Note that 
$G_V$ is the centralizer group $G_n(s)$ in this case. 
\end{exmp}

By linear algebra, the isomorphism classes of nilpotent 
representations are labelled by ($\mathbb{Z}/e\mathbb{Z}$-valued) multisegments. 

\begin{defn}
Let $l\in \Z_{>0}$ and $i\in\Z/e\Z$. The \emph{segment of
length $l$ and head $i$} is the sequence of consecutive residues 
$[i,i+1,...,i+l-1]$. We denote it by $[i;l)$. Similarly, 
The \emph{segment of length $l$ and tail $i$} is the sequence 
of consecutive residues 
$[i-l+1,...,i-1,i]$. We denote it by $(l;i]$. We say that 
$[i;l)$ has a \emph{(left) removable $i$-node} and 
$[i+1;l)$ has an \emph{(left) addable $i$-node}. 

A collection of segments is called a \emph{multisegment}. 
If the collection is the empty set, we call it the empty multisegment. 
\end{defn}

Each $[i;l)$ defines an indecomposable nilpotent $\C\Gamma_e$-module $\C[i;l)$, 
which is characterized by the property that 
\begin{quote}
$\C[i;l)$ is a uniserial module and $\Top(\C[i;l))=S_i$. 
\end{quote}
Hence, a complete set of isomorphism classes of nilpotent representations is 
given by the modules 
$$
M_\psi=\bigoplus_{i\in \Z/e\Z,l\in \Z_{>0}}\C[i,l)^{\oplus m_{[i; l)}},
$$
which is labelled by the multisegment
$$
\psi=\{[i; l)^{\oplus m_{[i; l)}}\}_{i\in \Z/e\Z,l\in \Z_{>0}}.
$$
We denote the corresponding $G_V$-orbit in $\Nil_V$ by $\mathcal O_\psi$. 

Now, we introduce the Hall polynomials. Let $\mathbb F_q$ be a finite field, 
and consider $\mathbb F_q\Gamma_e$-modules. Then, they are classified by 
multisegments again. 
Let $V$, $T$ and $W$ be $\Z/e\Z$-graded vector spaces over $\mathbb F_q$ such that 
$$
\underline{\dim}\;V=\underline{\dim}\;T+\underline{\dim}\;W.
$$
Let $\varphi_1$, $\varphi_2$ and $\psi$ be multisegments such that 
$\mathcal O_{\varphi_1}\subseteq \Nil_T$, 
$\mathcal O_{\varphi_2}\subseteq \Nil_W$ and 
$\mathcal O_\psi\subseteq \Nil_V$. If the number of submodules $U$ of $M_\psi$ 
that satisfies $U\simeq M_{\varphi_2}$ and $M_\psi/U\simeq M_{\varphi_1}$ is 
polynomial in $q=\mathrm{card}(\mathbb F_q)$, then this polynomial 
is called the \emph{Hall polynomial} and we denote it by 
$F_{\varphi_1,\varphi_2}^\psi(q)$. The existence of Hall polynomials 
in our case was proved by Jin Yun Guo \cite[Theorem 2.7]{Guo}. 

For $a$ and $b$ in $\Z^e$ we define a bilinear form $m$ by
$$
m(a,b)=\sum_{i\in \Z/e\Z}(a_ib_{i+1}+a_ib_i).
$$
We remark that this is not the Euler form used by Ringel to define 
his (twisted) Hall algebra, but the one 
used by Lusztig, which comes from the difference of dimensions of 
the fibers of two fiber bundles which appear in his geometric 
definition of the product, namely in the definition of 
the induction functor. In his theory, the Euler form appears 
in the definition of coproduct, namely in the definition of 
the restriction functor. 

Now, Lusztig's version of the \emph{Hall algebra} associated to 
$\Gamma_e$ is the 
$\C(v)$-algebra with basis $\{u_\psi \mid \text{$\psi$ is a multisegment}\}$ and
product is given by 
$$
u_{\varphi_1} u_{\varphi_2}=
v^{m(\underline\dim\;T,\;\underline\dim\;W)}
\sum_\psi F_{\varphi_1,\varphi_2}^\psi(v^{-2})u_\psi.
$$
Note that $[i;1)$ is the multisegment which labels the simple module $S_i$, for 
$i\in\Z/e\Z$. 
Then the $\C(v)$-subalgebra generated by
these $u_{[i;1)}$ is called the \emph{composition algebra}, 
and we may and do identify it with $U_v^-$ by $u_{[i;1)}\mapsto f_i$. 
For the proof, see \cite[Theorem 1.20]{Lu}. 

\begin{defn}
For each multisegment $\psi$, we define 
$E_\psi=v^{\dim \mathcal O_\psi}u_\psi$. The set 
$\{E_\psi \mid \text{$\psi$ is a multisegment.}\}$ 
is called the PBW basis of the Hall algebra. 
\end{defn}

\begin{exmp}
Let $e=3$. Then we have 
$$
f_1f_2=E_{\{[1;2)\}}+vE_{\{[1;1),[2;1)\}},\quad
f_2f_1=E_{\{[1;1),[2;1)\}}.
$$
Similarly, we have
\begin{equation*}
\begin{split}
f_0f_1f_2=&E_{\{[0;3)\}}+vE_{\{[1;2),[0;1)\}}
+vE_{\{[0;2),[2;1)\}}+v^2E_{\{[0;1),[1;1),[2;1)\}},\\
f_2f_0f_1=&E_{\{[2;3)\}}+vE_{\{[0;2),[2;1)\}}
+vE_{\{[2;2),[1;1)\}}+v^2E_{\{[2;1),[0;1),[1;1)\}},\\
f_1f_2f_0=&E_{\{[1;3)\}}+vE_{\{[2;2),[1;1)\}}
+vE_{\{[1;2),[0;1)\}}+v^2E_{\{[1;1),[2;1),[0;1)\}},\\
f_2f_1f_0=&E_{\{[2;2),[1;1)\}}+vE_{\{[0;1),[1;1),[2;1)\}},\\
f_0f_2f_1=&E_{\{[0;2),[2;1)\}}+vE_{\{[0;1),[1;1),[2;1)\}},\\
f_1f_0f_2=&E_{\{[1;2),[0;1)\}}+vE_{\{[0;1),[1;1),[2;1)\}}.
\end{split}
\end{equation*}
\end{exmp}
Note that $E_{\{[0;1),[1;1),[2;1)\}}$ does not appear with coefficient $1$. 
This is general phenomenon and we need aperiodicity to describe it. 

\begin{defn}
A multisegment $\psi$ is \emph{aperiodic} if, for every $l\in\Z_{>0}$, 
there exists some $i\in\Z/e\Z$ such that 
the segment of length $l$ and head $i$ does not appear 
in $\psi$. Equivalently, a multisegment $\psi$ is aperiodic if, 
for each $l\in\Z_{>0}$, there exists some $i\in\Z/e\Z$ such that 
the segment of length $l$ and tail $i$ does not appear 
in $\psi$. 
\end{defn}

The notion of aperiodicity and the following theorem are due to Lusztig. 
See \cite[15.3]{Lus} and \cite[Theorem 5.9]{L1}. 

\begin{thm}
\label{aperiodicity}
For each $b\in B(\infty)$, the canonical basis element $G_v(b)$ has the form 
$$
G_v(b)=E_\psi+\sum_{\psi'\ne \psi}c_{\psi,\psi'}(v)E_{\psi'},
$$
for a unique aperiodic multisegment $\psi$, such that 
$c_{\psi,\psi'}(v)\in\C(v)$ is regular at $v=0$ and $c_{\psi,\psi'}(0)=0$. 
\end{thm}

Hence, we may label elements of $B(\infty)$ by aperiodic 
multisegments. We identify $B(\infty)$ with the set of 
aperiodic multisegments. Then, we denote the canonical basis by 
$G_v(\psi)$, for multisegments $\psi$, hereafter. 

Leclerc, Thibon and Vasserot described the crystal structure 
on the set of aperiodic multisegments $B(\infty)$ 
in \cite[Theorem 4.1]{LTV}, by using a result by Reineke. 

Let $\psi$ be a multisegment. Let 
$\psi_{\ge l}$ be the multisegment obtained 
from $\psi$ by deleting multisegments of length less than $l$, 
for $l\in\Z_{>0}$. 
Let $m_{[i;l)}$ be the multiplicity of 
$[i;l)$ in $\psi$. Then, for $i\in \mathbb{Z}/e\mathbb{Z}$, 
we consider  
\begin{equation*}
S_{l,i}=\sum_{k\geq l}(m_{[i+1;k)}-m_{[i;k)}),
\end{equation*}
that is, the number of addable $i$-nodes of $\psi_{\ge l}$ 
minus the number of removable $i$-nodes of $\psi_{\ge l}$. 
Let $\ell_0<\ell_1<\cdots$ be those $l$ that attain 
$\min_{l>0}S_{l,i}$. The following is the description of 
the crystal structure given by Leclerc, Thibon and Vasserot. 

\begin{thm}
\label{LTV-1}
Let $\psi$ be a multisegment, $i\in\Z/e\Z$ and let $\ell_0$ be as above. Then, 
$\tilde f_i\psi=\psi_{\ell_0,i}$, where 
$\psi_{\ell_0,i}$ is obtained from $\psi$ by adding $[i;1)$ if $\ell_0=1$, 
and by replacing $[i+1;\ell_0-1)$ with $[i;\ell_0)$ if $\ell_0>1$. 
\end{thm}

\subsection{An anti-automorphism of $U_v^-$}

As the identification of the affine Hecke algebra with the convolution 
algebra $K^{G_n\times\C^\times}(Z_n)$, which will be explained in the 
next section, is not canonical, we go back and forth between 
two identifications. 
For this reason, we need another labelling by aperiodic multisegments. 

Let $V=\oplus_{i\in\Z/e\Z} V_i$ be a graded vector space as before, and 
define its dual graded vector space by 
$V^*=\oplus_{i\in\Z/e\Z} V^*_i$ where $V^*_i=\Hom_\C(V_{-i},\C)$. 
Then, by sending $X\in E_V$ to its transpose, we have a linear isomorphism
$$
\rho: E_V\simeq E_{V^*}=\oplus_{i\in\Z/e\Z}\Hom_\C(V^*_i,V^*_{i+1}).
$$
Using the standard basis of $E_V$ and its dual basis 
in $E_{V^*}$, we identify the underlying spaces $E_V$ and $E_{V^*}$. 
Note that the $G_V$-action on this $E_V$ is the conjugation by 
the transpose inverse of $g\in G_V$, while the $G_V$-action on 
the original $E_V$ is 
the conjugation by $g\in G_V$. Then, $\rho$ is an isomorphism of two 
$G_V$-varieties $E_V$ so that the $G_V$-orbit $\mathcal O_\psi$ 
in the original $E_V$ corresponds to the $G_V$-orbit 
$\mathcal O_{\rho(\psi)}$ in the new $E_V$, where $\rho(\psi)$ is 
defined by $\rho([i;l))=(l;-i]$. Thus, 
we have a linear isomorphism of the Hall algebras on both sides, 
which we also denote by $\rho$, such that
$$
\rho(E_\psi)=E_{\rho(\psi)}\;\;\text{and}\;\;
\rho(G_v(\psi))=G_v(\rho(\psi))\;\text{if $\psi$ is aperiodic.}
$$
That is, this gives a relabelling 
of the PBW basis and the canonical basis. However, 
if we take the algebra structure into account, $\rho$ induces 
the anti-automorphism of $U_v^-$ given by $f_i\mapsto f_{-i}$, which 
is clear from the definition of the multiplication of the Hall algebra. 
In particular, 
the crystal structure on the set of aperiodic multisegments 
is changed in this new labelling, and the Kashiwara 
operators $\tilde e_i$ and $\tilde f_i$ correspond to 
the Kashiwara operators 
$\tilde e_{-i}$ and $\tilde f_{-i}$ in this new crystal structure. 
In the new crystal structure, 
we change the definition of addable and removable $i$-nodes as follows. 

\begin{defn}
We say that $(l;i]$ has a \emph{(right) removable $i$-node} and 
$(l;i-1]$ has an \emph{(right) addable $i$-node}. 
\end{defn}

We consider  
$S_{l,i}=\sum_{k\geq l}(m_{(k;i-1]}-m_{(k;i]})$, 
that is, the number of addable $i$-nodes of $\psi_{\ge l}$ 
minus the number of removable $i$-nodes of $\psi_{\ge l}$ 
in the new definition of removable and addable $i$-nodes. 
Let $\ell_0<\ell_1<\cdots$ be those $l$ that attain 
$\min_{l>0}S_{l,i}$. Then, the crystal structure in the new 
labelling is given as follows. In fact, this version is stated in \cite{LTV}. 

\begin{thm}
\label{LTV-2}
Let $\psi$ be a multisegment, $i\in\Z/e\Z$ and let $\ell_0$ be as above. Then, 
$\tilde f_i\psi=\psi_{\ell_0,i}$, where 
$\psi_{\ell_0,i}$ is obtained from $\psi$ by adding $(1;i]$ if $\ell_0=1$, 
and by replacing $(\ell_0-1;i-1]$ with $(\ell_0;i]$ if $\ell_0>1$. 
\end{thm}

To compute $\tilde e_i\psi$, for a multisegment $\psi$, 
we consider the same $S_{l,i}$. If 
$\min_{l>0} S_{l,i}=0$, 
then $\tilde e_i\psi=0$. Otherwise, let $\ell_0$ be the maximal $l$ 
that attains $\min_{l>0} S_{l,i}$. Then, $\tilde e_i\psi$ is obtained 
from $\psi$ by replacing $(\ell_0;i]$ with $(\ell_0-1;i-1]$. 

We use the crystal structure on the set of aperiodic multisegments in 
Theorem \ref{LTV-1} when we choose the identification of 
$R(G_n\times\C^\times)$-algebras 
$\Hecke_n\simeq K^{G_n\times\C^\times}(Z_n)$ following Lusztig \cite{L3}, 
while we use that in Theorem \ref{LTV-2} when we choose 
the identification 
$\Hecke_n\simeq K^{G_n\times\C^\times}(Z_n)$ following Ginzburg \cite{CG}. 
We note that the second crystal structure is the star crystal 
structure of the first.

\section{Affine Hecke algebras}

Let $\Hecke_n$ be the extended affine Hecke algebra associated with $G_n$. It is the $\C[q^{\pm1}]$-algebra generated by 
$T_i$, for $1\le i<n$, and $X_i^{\pm1}$, for $1\le i\le n$, subject to the relations
$$
(T_i-q)(T_i+1)=0,\quad q^{-1}T_iX_iT_i=X_{i+1},\quad\text{etc.}
$$
In this section, we recall the geometric realization of affine Hecke algebras 
by Lusztig and Ginzburg, and of specialized affine Hecke algebras by Ginzburg. 

\subsection{Varieties}
Let $G_n=GL_n(\C)$ as before, and $B_n$ the Borel subgroup of upper 
triangular matrices. We denote the unipotent radical of $B_n$ by $U_n$, and
the maximal torus of diagonal matrices by $T_n$. 
Write $\C^n=\C e_1\oplus\cdots\oplus\C e_n$ and let 
$\Fl_n$ be the flag variety, which consists of increasing subspaces 
$F=(F_i)_{0\le i\le n}$ in $\C^n$ such that $\dim F_i=i$, for all $i$.
We consider the diagonal $G_n$-action on $\Fl_n\times\Fl_n$. Then, $G_n$-orbits in $\Fl_n\times\Fl_n$ 
are in bijection with $B_n$-orbits in $\Fl_n$ and if we denote
$$
\{(F,F')\in \Fl_n\times\Fl_n\mid 
\dim(F_i\cap F'_j)=\sharp\{k \mid 1\le k\le i,\;1\le w(k)\le j\}\}
$$
by $O_n(w)$, for $w\in \Sym_n$, they give a complete set of $G_n$-orbits, and 
a pair of flags $(F,F')$ belongs to $O_n(w)$ if and only if
$$
\dim\frac{F_i\cap F'_j}{F_{i-1}\cap F'_j+F_i\cap F'_{j-1}}=
\begin{cases} 1 \;\;&(j=w(i))\\ 0 \;\;&(\text{otherwise})\end{cases}.
$$
We denote by $\Nil_n$ the set of nilpotent elements in $Mat_n(\C)$ and 
write
$$
Y_n=\{(X,F)\in \Nil_n\times\Fl_n \mid XF_i\subseteq F_{i-1}\}\simeq T^*\Fl_n. 
$$
Then the Steinberg variety is defined by
\begin{equation*}
\begin{split}
Z_n&=Y_n\times_{\Nil_n}Y_n\\
&=\{(X,F,F')\in \Nil_n\times\Fl_n\times\Fl_n 
\mid XF_i\subseteq F_{i-1}, XF'_i\subseteq F'_{i-1}\}.
\end{split}
\end{equation*}
$Z_n$ is a $G_n\times\C^\times$-variety by the action
$$
(g,c)(X,F,F')=(c^{-1}gXg^{-1},gF,gF'),
$$
for $(g,c)\in G_n\times\C^\times$ and $(X,F,F')\in Z_n$. 

We consider the complexified K-group of the abelian category of 
$G_n\times\C^\times$-equivariant coherent sheaves on $Z_n$. 
Using the closed embedding $Z_n\subseteq Y_n\times Y_n$, 
we have the convolution algebra 
$K^{G_n\times\C^\times}(Z_n)$. $Z_n$ has a partition 
$Z_n=\sqcup_{w\in\Sym_n}Z_n(w)$, where
$$
Z_n(w)=\{(X,F,F')\in Z_n \mid (F,F')\in O_n(w)\}.
$$
We have $\dim Z_n(w)=n(n-1)$ and $Z_n(w)$ is a 
$(\frac{n(n-1)}{2}-\ell(w))$-dimensional vector bundle over $O_n(w)$. 
Then, $\{\overline{Z_n(w)}\}_{w\in\Sym_n}$ is 
the set of the irreducible components of $Z_n$. Define
$$
Z_{n-1,n}=\{(X,F,F')\in Z_n \mid F_{n-1}=F'_{n-1}\}.
$$
The condition $F_{n-1}=F'_{n-1}$ is equivalent to 
$(F,F')\in\sqcup_{w\in\Sym_{n-1}}O_n(w)$, so that we have 
$Z_{n-1,n}=\sqcup_{w\in\Sym_{n-1}}Z_n(w)$. 

Similarly, $(F,F')\in O_n(e)\sqcup O_n(s_i)=\overline{O_n(s_i)}$ 
if and only if $F_j=F'_j$, for all $j\ne i$, and 
$$
\overline{Z_n(s_i)}=\{(X,F,F')\in Z_n \mid F_j=F'_j,\;\text{for all $j\ne i$},\;
XF_{i+1}\subseteq F_{i-1}\}.
$$
The pushforward of $\mathcal O_{\overline{Z_n(s_i)}}$ 
with respect to the closed embedding $\overline{Z_n(s_i)}\subseteq Z_n$ is also 
denoted by $\mathcal O_{\overline{Z_n(s_i)}}$ by abuse of notation. We denote
$$
b_i=[\mathcal O_{\overline{Z_n(s_i)}}]\in K^{G_n\times\C^\times}(Z_n).
$$

Let $Q_{i,i+1}$ be the parabolic subgroup of $G_n$ which corresponds to $s_i$, 
$\mathfrak n_{i,i+1}$ the nilradical of its Lie algebra. Then
$$
\overline{Z_n(s_i)}=(G_n\times\C^\times)\times_{Q_{i,i+1}\times\C^\times}
(\mathfrak n_{i,i+1}\times \mathbb P^1\times\mathbb P^1)
$$
is a vector bundle over 
$\overline{O_n(s_i)}=(G_n\times\C^\times)\times_{Q_{i,i+1}\times\C^\times}
(\mathbb P^1\times\mathbb P^1)$. Then we define as follows. 

\begin{defn}
The line bundle $\mathcal L_i$ on $\overline{Z_n(s_i)}$ is the pullback of 
$$
(G_n\times\C^\times)\times_{Q_{i,i+1}\times\C^\times}
(\mathcal O_{\mathbb P^1}(-1)\otimes\mathcal O_{\mathbb P^1}(-1))
$$
on $\overline{O_n(s_i)}$. 
\end{defn}

For $\lambda\in \Z\epsilon_1\oplus\cdots\oplus\Z\epsilon_n=\Hom(T_n,\C^\times)$, let $\C_\lambda$ be the $B_n\times\C^\times$-module associated with $\lambda$ and define
the associated line bundle $L_\lambda$ on $\Fl_n$ by
$$
L_\lambda=(G_n\times\C^\times)\times_{B_n\times\C^\times}\C_\lambda.
$$
When we consider $\lambda$ as a character of $T_n$, we denote it 
by $e^\lambda$. Then, we identify 
$K^{G_n\times\C^\times}(\Fl_n)=R(T_n\times\C^\times)$ via 
$L_\lambda\mapsto e^\lambda$ as usual. 

Let us denote $\pi_n:Y_n\rightarrow \Fl_n$ and 
$\delta_n: Z_n(e)\subseteq Z_n$. We consider the diagram
$$
\Fl_n\overset{\pi_n}{\longleftarrow} Y_n\simeq Z_n(e)
\overset{\delta_n}{\longrightarrow} Z_n
$$
and we denote
$$
\theta_\lambda=[{\delta_n}_*\pi_n^*L_{-\lambda}]\in K^{G_n\times\C^\times}(Z_n). 
$$

\begin{defn}
We define $T_i=[\mathcal L_i]+q$, for $1\le i<n$, 
and $X_i=\theta_{\epsilon_i}$, for $1\le i\le n$. 
\end{defn}

We have $\theta_\lambda=\prod_{i=1}^n X_i^{\lambda_i}$, for 
$\lambda=\sum_{i=1}^n\lambda_i\epsilon_i$. Using the exact sequence
$$
0\rightarrow \mathcal O_{\mathbb P^1}(-1)\otimes\mathcal O_{\mathbb P^1}(-1)
\rightarrow \mathcal O_{\mathbb P^1}\otimes\mathcal O_{\mathbb P^1}
\rightarrow \mathcal O_{\Delta\mathbb P^1}\rightarrow 0
$$
where $\Delta\mathbb P^1\subseteq \mathbb P^1\times\mathbb P^1$ is the diagonal, 
we know that $[\mathcal L_i]=b_i-(1-q\theta_{\alpha_i})$. 

Then, $T_i$, for $1\le i<n$, and $X_i^{\pm1}$, for $1\le i\le n$, satisfy the defining relations of $\Hecke_n$. In particular, we have the Bernstein relation
$$
T_i\theta_\lambda=\theta_{s_i\lambda}T_i+(1-q)
\frac{\theta_\lambda-\theta_{s_i\lambda}}{\theta_{-\alpha_i}-1},
$$
where $\alpha_i=-\epsilon_i+\epsilon_{i+1}$. 
This follows from the next theorem. The theorem was found by Lusztig and 
the action of $T_i$ is called the Demazure-Lusztig operator. 

\begin{thm}
Through the Thom isomorphism, we identify $K^{G_n\times\C^\times}(Y_n)$ with 
$$
K^{G_n\times\C^\times}(\Fl_n)= R(T_n\times\C^\times).
$$
Then the convolution action of $K^{G_n\times\C^\times}(Z_n)$ on 
$K^{G_n\times\C^\times}(Y_n)$ is given by
$$
T_if=\frac{f-s_if}{e^{\alpha_i}-1}-q\frac{f-e^{\alpha_i}s_if}{e^{\alpha_i}-1},\;\;
X_if=e^{-\epsilon_i}f.
$$
\end{thm}

It is well-known that this is a faithful representation of $\Hecke_n$. Note that 
we have chosen the isomorphism 
$\Hecke_n\simeq K^{G_n\times\C^\times}(Z_n)$ to have  
the same formulas as \cite[Theorem 7.2.16, Proposition 7.6.38]{CG}. 
When we follow \cite{L3}, we define
$$
\theta_\lambda=[{\delta_n}_*\pi_n^*L_\lambda]\;\;\text{and}\;\;
T_i=-[\mathcal L_i]-1. 
$$
Then, the formulas for the convolution action 
on $R(T_n\times\C^\times)$ change to those in 
\cite[p.335]{L3}. The two identifications of 
$\Hecke_n\simeq K^{G_n\times\C^\times}(Z_n)$ are related by the involution $\sigma$ 
defined by 
$$
T_i\mapsto -qT_i^{-1}, \;\;X_i\mapsto X_i^{-1}.
$$
\textit{In the rest of this section, we follow the identification in \cite{L3}}. 


The center $Z(\Hecke_n)$ of $\Hecke_n$ is the $\C[q^{\pm1}]$-subalgebra consisting 
of all the symmetric Laurent polynomials in $X_1,\dots,X_n$. 
Thus, we identify $Z(\Hecke_n)$ with $R(G_n\times\C^\times)$. We also identify 
$\C[q^{\pm1}][X_1^{\pm1},\dots,X_n^{\pm1}]$ with $R(T_n\times\C^\times)$. 

Let $K^{G_n\times\C^\times}(Z_{n-1,n})$ be the convolution algebra with respect to the embedding $Z_{n-1,n}\subseteq Y_n\times Y_n$. Let
\begin{itemize}
\item
$\Hecke_{n-1,n}$ be the parabolic subalgebra $\Hecke_{n-1}\otimes_\C\C[X_n^{\pm1}]$ of $\Hecke_n$, and
\item
$\iota_n:Z_{n-1,n}\subseteq Z_n$ be the inclusion map. 
\end{itemize}

We attribute the next theorem to Ginzburg \cite{CG} and Lusztig \cite{L3}. 
In \cite{KL}, it was stated as an isomorphism of bimodules. 

\begin{thm}
\label{geometric realization}
\item[(1)]
We have an isomorphism of $R(G_n\times\C^\times)$-algebras
$\Hecke_n\simeq K^{G_n\times\C^\times}(Z_n)$ by the above choice of $T_i$ and $X_i$ 
in $K^{G_n\times\C^\times}(Z_n)$. 
\item[(2)]
The inclusion map $\iota_n$ induces 
the following commutative diagram of $Z(\Hecke_n)$-algebras. 
\begin{eqnarray*}
{\iota_n}_*:K^{G_n\times\C^\times}(Z_{n-1,n})&\rightarrow& K^{G_n\times\C^\times}(Z_n)\\
\downarrow\qquad & & \quad\downarrow \\
\Hecke_{n-1,n}&\subseteq& \;\;\Hecke_n
\end{eqnarray*}
where the vertical arrows are isomorphisms. 
\end{thm}

It is also clear that the inclusion map 
$Y_n\simeq Z_n(e)\hookrightarrow Z_{n-1,n}$ induces 
$$
K^{G_n\times\C^\times}(Y_n)\rightarrow K^{G_n\times\C^\times}(Z_{n-1,n})
$$
and it is identified with 
$R(T_n\times\C^\times) \hookrightarrow \Hecke_{n-1,n}$. 

\subsection{The embedding of $H_{n-1}$ into $H_n$}
Let 
$$
Y_{n-1,n}=T^*\Fl_n|_{\Fl_{n-1}}, 
$$
where we identify 
$\Fl_{n-1}=\{ F\in\Fl_n \mid F_{n-1}=\C^{n-1}\}$, and let 
$$
\Nil_{n-1,n}=\{X\in\Nil_n \mid X\C^{n-1}\subseteq\C^{n-1}\}.
$$
Then we define 
\begin{equation*}
\begin{split}
Z'_{n-1,n}&=Y_{n-1,n}\times_{\Nil_{n-1,n}}Y_{n-1,n}\\
&=\{(X,F,F')\in Z_n \mid F_{n-1}=F'_{n-1}=\C^{n-1}\}.
\end{split}
\end{equation*}

Let $P_{n-1,n}$ be the maximal parabolic subgroup of $G_n$ that stabilizes $\C^{n-1}$. 
The Levi part $L_{n-1,n}\times\C^\times$ of $P_{n-1,n}\times\C^\times$ is 
$(G_{n-1}\times\C^\times)\times\C^\times$, 
which acts on $Z_{n-1}$ by letting the middle component act trivially. We denote the unipotent radical of $P_{n-1,n}$ by $U_{n-1,n}$. It is also the unipotent radical 
of $P_{n-1,n}\times\C^\times$. Explicitly, 
$$
L_{n-1,n}=\begin{pmatrix} G_{n-1} & 0 \\ 0 & \C^\times\end{pmatrix},\quad
U_{n-1,n}=\left\{\begin{pmatrix} 1_{n-1} & * \\ 0 & 1\end{pmatrix}\right\}. 
$$
We consider the following diagram.
$$
Z'_{n-1,n}\overset{\mu_{n-1,n}}{\longleftarrow}
(G_n\times\C^\times)\times Z'_{n-1,n}
\overset{\nu_{n-1,n}}{\longrightarrow}
(G_n\times\C^\times)\times_{P_{n-1,n}\times\C^\times}Z'_{n-1,n}=Z_{n-1,n}.
$$
Then we have the restriction map
$$
\Res^{G_n\times\C^\times}_{P_{n-1,n}\times\C^\times}:
K^{G_n\times\C^\times}(Z_{n-1,n})\simeq K^{P_{n-1,n}\times\C^\times}(Z'_{n-1,n}).
\label{eq1}
$$
$Z'_{n-1,n}$ is a $L_{n-1,n}\times\C^\times$-equivariant vector bundle of rank $n-1$
over $Z_{n-1}$ and we write $\kappa_{n-1,n}:Z'_{n-1,n}\rightarrow Z_{n-1}$. Then
$\kappa_{n-1,n}^*$ gives the Thom isomorphism
$$
K^{L_{n-1,n}\times\C^\times}(Z'_{n-1,n})\overset{\sim}{\leftarrow}
K^{L_{n-1,n}\times\C^\times}(Z_{n-1}).
$$
Noting that 
$$
K^{P_{n-1,n}\times\C^\times}(Z'_{n-1,n})
\simeq K^{L_{n-1,n}\times\C^\times}(Z'_{n-1,n})
$$
by the forgetful map, and
$$
K^{L_{n-1,n}\times\C^\times}(Z_{n-1})
\simeq K^{G_{n-1}\times\C^\times}(Z_{n-1})\otimes_\C\C[X_n^{\pm1}],
$$
we have
$$
K^{P_{n-1,n}\times\C^\times}(Z'_{n-1,n})\simeq K^{G_{n-1}\times\C^\times}(Z_{n-1})\otimes_\C\C[X_n^{\pm1}].
$$
Now, the following holds. 

\begin{prop}
\label{the first isom}
We have the following 
isomorphism of $R(L_{n-1,n}\times\C^\times)$-algebras
$$
K^{G_n\times\C^\times}(Z_n)\supseteq K^{G_n\times\C^\times}(Z_{n-1,n})
\simeq K^{G_{n-1}\times\C^\times}(Z_{n-1})\otimes_\C\C[X_n^{\pm1}],
$$
which gets identified with 
$\Hecke_n\supseteq\Hecke_{n-1,n}=\Hecke_{n-1}\otimes_\C\C[X_n^{\pm1}]$. 
\end{prop}
\begin{proof}
We only have to show that $b_i\mapsto b_i$, for $1\le i<n-1$, and 
$\theta_\lambda\mapsto \theta_\lambda$. Define 
$$
\overline{Z'_{n-1,n}(s_i)}=\{(X,F,F')\in Z'_{n-1,n}\mid 
F_j=F'_j\;\text{for all $j\ne i$,}\;\;XF_{i+1}\subseteq F_{i-1}\}.
$$
Then $\nu_{n-1,n}^{-1}(\overline{Z_n(s_i)})
=(G_n\times\C^\times)\times \overline{Z'_{n-1,n}(s_i)}$ and we have
$$
\nu_{n-1,n}^*\mathcal O_{\overline{Z_n(s_i)}}=
\mu_{n-1,n}^*\mathcal O_{\overline{Z'_{n-1,n}(s_i)}},\;\;
\mathcal O_{\overline{Z'_{n-1,n}(s_i)}}=
\kappa_{n-1,n}^*\mathcal O_{\overline{Z_{n-1}(s_i)}}.
$$
Hence, $b_i\mapsto b_i$, for $1\le i<n-1$.

Let $Z'_{n-1,n}(e)=\{(X,F,F')\in Z'_{n-1,n}\mid F=F'\}$ and consider the diagram
\begin{eqnarray*}
& & \nu_{n-1,n}^{-1}(Z_n(e))=(G_n\times\C^\times)\times Z'_{n-1,n}(e)\\
& \swarrow     & \qquad\downarrow \nu_{n-1,n}\\
Z'_{n-1,n}(e)\simeq Y_{n-1,n} &\;\subseteq & \qquad Y_n\simeq Z_n(e)\\ 
\pi_{n-1,n}\downarrow\quad       &          & \quad\quad\downarrow\pi_n \\
 \Fl_{n-1}            &\;\subseteq & \quad\quad\Fl_n
\end{eqnarray*}
Then $\nu_{n-1,n}^*\pi_n^*L_\lambda=\mu_{n-1,n}^*\pi_{n-1,n}^*L_\lambda|_{\Fl_{n-1}}$
and 
$$
L_\lambda|_{\Fl_{n-1}}=(P_{n-1,n}\times\C^\times)\times_{B_n\times\C^\times}\C_\lambda.
$$
But the diagram
\begin{eqnarray*}
Z'_{n-1,n}(e)=\kappa_{n-1,n}^{-1}(Z_{n-1}(e))\simeq Y_{n-1,n}& 
\overset{\kappa_{n-1,n}}{\longrightarrow} & Y_{n-1}\simeq Z_{n-1}(e) \\
&\searrow & \;\downarrow\pi_{n-1} \\
&         & \Fl_{n-1}             
\end{eqnarray*}
shows
$$
\pi_{n-1,n}^*L_\lambda|_{\Fl_{n-1}}=
\kappa_{n-1,n}^*\left(\pi_{n-1}^*((P_{n-1,n}\times\C^\times)\times_{B_n\times\C^\times}\C_\lambda)\right).
$$
Hence, $\theta_\lambda\mapsto \theta_\lambda$, for $\lambda\in\Hom(T_n,\C^\times)$. 

As the generators $b_i$ and $\theta_\lambda$ correspond correctly, it is an 
isomorphism of $R(L_{n-1,n}\times\C^\times)$-algebras, which is identified with  
$\Hecke_{n-1}\otimes\C[X_n^{\pm1}]\hookrightarrow \Hecke_n$.
\end{proof}

\subsection{Specialized Hecke algebras}

Let $\zeta\in\C$ be a primitive $e^{th}$ root of unity, for $e\ge2$. 
We fix a diagonal matrix $s=\diag(\zeta^{s_1},\dots,\zeta^{s_n})$, and set
$a=(s,\zeta)\in G_n\times\C^\times$. 
We denote by $A$ the smallest closed algebraic subgroup of 
$G_n\times\C^\times$ that contains $a$, namely the cyclic group $\langle a\rangle$
of order $e$ in our case. Note that $A$ is contained in 
$(G_{n-1}\times\C^\times)\times\C^\times$. 

\begin{defn}
We denote the $A$-fixed points of $M$ by $M^a$, for $M=Z_n$, $Z_{n-1,n}$, 
$Y_n=T^*\Fl_n$, $\Fl_n$ etc. 
\end{defn}

Let $\C_a$ be the $R(T_n\times\C^\times)$-module defined by 
$X_i\mapsto \zeta^{s_i}$, for $1\le i\le n$, and $q\mapsto \zeta$. 
$\C_a|_{R(G_n\times\C^\times)}$ defines a central character
$Z(H_n)=R(G_n\times\C^\times)\rightarrow \C$. Then we write
$\C_a\otimes_{Z(H_n)}-$, for the specialization of 
the center with respect to the central character. We define $f_a\in \C[X_n]$ by
$$
f_a(X_n)=(X_n-\zeta^{s_1})\cdots(X_n-\zeta^{s_n}).
$$

\begin{defn}
The $\C$-algebra $\Hecke_n^a=\C_a\otimes_{Z(\Hecke_n)}\Hecke_n$ is called the 
\emph{specialized Hecke algebra} of rank $n$ at $a$. 
The specialized algebra $\C_a\otimes_{Z(\Hecke_n)}\Hecke_{n-1,n}$ 
of the parabolic subalgebra $\Hecke_{n-1,n}$ is denoted 
$\Hecke_{n-1,n}^a$.
\end{defn} 

\begin{lem}
\label{structure lemma}
Let $a_k$ be the $k$-th elementary symmetric function in $X_1,\dots,X_n$ evaluated at $X_1=\zeta^{s_1},\dots,X_n=\zeta^{s_n}$.
Then, $\C[X_n^{\pm1}]/(f_a)$ is a $Z(H_{n-1})$-algebra via
$$
e_k\mapsto a_k-a_{k-1}X_n+\cdots+(-1)^kX_n^k,
$$
where $e_k$ is the $k$-th elementary symmetric function in 
$X_1,\dots, X_{n-1}$, and 
$$
H_{n-1,n}^a=H_{n-1}\otimes_{Z(H_{n-1})}\C[X_n^{\pm1}]/(f_a).
$$
\end{lem}
\begin{proof}
As $e_k+X_ne_{k-1}$ is the $k$-th elementary symmetric function in $X_1,\dots, X_n$, 
the surjective map
$$
H_{n-1,n}=H_{n-1}\otimes_\C\C[X_n^{\pm1}]\rightarrow 
H_{n-1}\otimes_{Z(H_{n-1})}\C[X_n^{\pm1}]/(f_a)
$$
factors through $H_{n-1,n}^a$. On the other hand, both have the same dimension 
$n!(n-1)!$. Hence the result. 
\end{proof}

As $\Hecke_n\simeq K^{G_n\times\C^\times}(Z_n)$ and 
$\Hecke_{n-1,n}\simeq K^{G_n\times\C^\times}(Z_{n-1,n})$ as 
$R(G_n\times\C^\times)$-algebras, 
we identify the following $\C$-algebras respectively.
\begin{align*}
\Hecke_n^a&=\C_a\otimes_{R(G_n\times\C^\times)}K^{G_n\times\C^\times}(Z_n),\\
\Hecke_{n-1,n}^a&= 
\C_a\otimes_{R(G_n\times\C^\times)}K^{G_n\times\C^\times}(Z_{n-1,n}).
\end{align*}

By Proposition \ref{the first isom}, 
geometric realization of Lemma \ref{structure lemma} is given by 
\begin{align*}
\C_a\otimes_{Z(H_n)}K^{G_n\times\C^\times}(Z_{n-1,n})
&\simeq \C_a\otimes_{Z(H_n)}
\left(K^{G_{n-1}\times\C^\times}(Z_{n-1})\otimes \C[X_n^{\pm1}]\right)\\
&\simeq K^{G_{n-1}\times\C^\times}(Z_{n-1})\otimes_{Z(H_{n-1})} \C[X_n^{\pm1}]/(f_a).
\end{align*}

Let $m_i$ be the multiplicity of $\zeta^i$ in 
$\{\zeta^{s_1},\dots,\zeta^{s_n}\}$. Then
$$
\C[X_n^{\pm1}]/(f_a)\simeq
\bigoplus_{i\in\Z/e\Z}\C[X_n^{\pm1}]/((X_n-\zeta^i)^{m_i}).
$$

\begin{defn}
We denote by $p_i$ the identity of $\C[X_n^{\pm1}]/((X_n-\zeta^i)^{m_i})$ which is 
viewed as an element of $\Hecke_{n-1,n}^a$. Thus, $p_i$ are central idempotents 
of $\Hecke_{n-1,n}^a$ such that $\sum_{i\in \Z/e\Z} p_i=1$ and 
$p_ip_j=p_jp_i=\delta_{ij}p_i$. 
\end{defn}

We have the decomposition of $\C$-algebras 
$$
\Hecke_{n-1,n}^a=\bigoplus_{i\in\Z/e\Z} p_i\Hecke_{n-1,n}^ap_i.
$$

We fix $i\in\Z/e\Z$ and suppose that $(\nu_1,\dots,\nu_n)$ is a permutation of 
$(s_1,\dots,s_n)$ such that $\nu_n=i$. Then 
$$
(\diag(\zeta^{\nu_1},\dots,\zeta^{\nu_{n-1}}),\zeta)\in G_{n-1}\times\C^\times
$$
defines a central character of $H_{n-1}$ and we may define the specialized 
Hecke algebra with respect to the central character. We denote it by 
$H_{n-1}^{a;i}$. 
By Lemma \ref{structure lemma}, we have the surjective algebra homomorphism
$$
p_i\Hecke_{n-1,n}^ap_i\longrightarrow \Hecke_{n-1}^{a;i},
$$
because if we write $b_k=a_k-a_{k-1}\zeta^i+\cdots+(-1)^k\zeta^{ki}$, then 
\begin{align*}
\left(\sum_{k=0}^{n-1} (-1)^kb_kT^k\right)(1-\zeta^{\nu_n}T)
&=\sum_{k=0}^n (-1)^k(b_{k-1}\zeta^i+b_k)T^k\\
&=\sum_{k=0}^n (-1)^ka_kT^k\\
&=(1-\zeta^{s_1}T)\cdots(1-\zeta^{s_n}T)
\end{align*}
and $e_k\mapsto b_k$, for $1\leq k\leq n-1$, is the central character 
which defines $H_{n-1}^{a;i}$. Composing it with the projection 
$\Hecke_{n-1,n}^a$ to $p_i\Hecke_{n-1,n}^ap_i$, we have
$$
\Hecke_{n-1,n}^a\longrightarrow \Hecke_{n-1}^{a;i},
$$
which is nothing but the specialization map at $X_n=\zeta^i$. Its 
geometric realization is given by 
$$
\C_a\otimes_{Z(H_n)}K^{G_n\times\C^\times}(Z_{n-1,n})
\longrightarrow 
\C_{a;i}\otimes_{Z(H_{n-1})}K^{G_{n-1}\times\C^\times}(Z_{n-1}).
$$

\begin{lem}
\label{simple modules}
Simple $p_i\Hecke_{n-1,n}^ap_i$-modules are obtained from simple 
$\Hecke_{n-1}^{a;i}$-modules through the algebra homomorphism 
$p_i\Hecke_{n-1,n}^ap_i\rightarrow \Hecke_{n-1}^{a;i}$.
\end{lem}
\begin{proof}
Let $I_a$ be the two-sided ideal of $p_i\Hecke_{n-1,n}^ap_i$ generated by 
$X_n-\zeta^i$. Then $I_a$ is nilpotent, so that $I_a$ acts as zero 
on simple $p_i\Hecke_{n-1,n}^ap_i$-modules. As 
$H_{n-1}^{a;i}=p_i\Hecke_{n-1,n}^ap_i/I_a$ by Lemma \ref{structure lemma}, 
we have the result. 
\end{proof}

Ginzburg's theory tells us how to realize the specialized Hecke algebra 
in sheaf theory. Definitions of the maps in (1) are necessary in the 
proof of (2), so that they will be given in the proof. 

\begin{thm}
\label{Borel-Moore}
\item[(1)]
We may identify $\Hecke^a_n=H_*^{BM}(Z_n^a,\C)$ by
$$
\C_a\otimes_{R(A)}K^{G_n\times\C^\times}(Z_n)\simeq
\C_a\otimes_{R(A)}K^A(Z_n)\overset{res_n}{\simeq}
K(Z^a_n)\overset{RR_n}{\simeq}H_*^{BM}(Z_n^a,\C).
$$
We may identify $\Hecke^a_{n-1,n}=H_*^{BM}(Z_{n-1,n}^a,\C)$ in the same way. 
\item[(2)]
The following diagram of $\C$-algebras commutes. 
\begin{eqnarray*}
H_*^{BM}(Z_{n-1,n}^a,\C) &\overset{{\iota_n}_*}{\longrightarrow}
& H_*^{BM}(Z_n^a,\C)\\
\parallel\qquad\qquad  & & \quad\qquad\parallel \\
\Hecke_{n-1,n}^a\quad\quad & \hookrightarrow & \;\quad\quad\Hecke_n^a
\end{eqnarray*}
Similarly, $Y^a_n\rightarrow Z^a_{n-1,n}$ induces 
the following commutative diagram of $\C$-algebras.
\begin{eqnarray*}
H_*^{BM}(Y^a_n,\C) & \longrightarrow
& H_*^{BM}(Z^a_{n-1,n},\C)\\
\parallel\qquad\quad  & & \quad\qquad\parallel \\
\C_a\otimes_{Z(H_n)}R(T_n\times\C^\times) & \hookrightarrow & 
\;\quad\quad\Hecke_{n-1,n}^a
\end{eqnarray*}
\end{thm}
\begin{proof}
(1) is well-known. See \cite{Ari5} or \cite{CG}. We check the commutativity in (2). 
Let $Y_n^a\times Y_n^a\subseteq Y_n\times Y_n$ be the closed embedding. 
Then, we have
\begin{equation*}
\begin{split}
i^A_n:&\;\; Z_n^a=Z_n\cap(Y_n^a\times Y_n^a) \hookrightarrow Z_n,\\
i^A_{n-1,n}:&\;\; Z_{n-1,n}^a=Z_{n-1,n}\cap(Y_n^a\times Y_n^a) \hookrightarrow Z_{n-1,n}. 
\end{split}
\end{equation*}
We define the pullback
\begin{eqnarray*}
(i^A_n)^*:& K^A(Z_n)\longrightarrow K^A(Z_n^a), \\
(i^A_{n-1,n})^*:& K^A(Z_{n-1,n})\longrightarrow K^A(Z_{n-1,n}^a), 
\end{eqnarray*}
in terms of the embedding $Y_n^a\times Y_n^a\subseteq Y_n\times Y_n$. 
We have the linear $A$-action on each fiber of the normal bundle 
$T_{Y_n^a}Y_n$ and its decomposition into isotropic components leads to 
the decomposition of the normal bundle into the direct sum of  
vector bundles $N_i$, for $i\in\Z/e\Z$, over $Y_n^a$. We define $\lambda_n$ by 
$$
\lambda_n=\bigotimes_{i\in\Z/e\Z}
\left(\sum_{j\ge0} (-\zeta^i)^j\wedge^jN_i^\vee\right)\in K(Y_n^a).
$$
$res_n$ for $Z_n$ is defined by 
$$
res_n: K^A(Z_n)\overset{(i^A_n)^*}{\longrightarrow} 
\C_a\otimes_{R(A)}K^A(Z_n^a) \simeq 
K(Z_n^a)\overset{1\otimes\lambda_n^{-1}}{\longrightarrow}K(Z_n^a),
$$
and similarly for $Z_{n-1,n}$ and $Y_n$. Here, 
$1\otimes\lambda_n^{-1}\in K(Y_n^a\times Y_n^a)$ acts on $K(Z_n^a)$ 
by the multiplication. Then, the commutativity of the diagram
\begin{eqnarray*}
\C_a\otimes_{R(A)}K^A(Y_n) &\simeq & \;\; K(Y_n^a) \\
\downarrow\qquad & & \qquad \downarrow \\
\C_a\otimes_{R(A)}K^A(Z_{n-1,n}) &\simeq & K(Z_{n-1,n}^a) \\
{\iota_n}_*\downarrow\qquad & & \qquad \downarrow{\iota_n}_* \\
\C_a\otimes_{R(A)}K^A(Z_n) &\simeq & \;\; K(Z_n^a) \\
\end{eqnarray*}
follows from the statement below. 
\begin{quote}
Let $N\subseteq M$ be a closed embedding between smooth varieties, 
$Z$ a closed subvariety of $M$. Let $Z'=Z\cap N$ and denote
\begin{eqnarray*}
N & \overset{\psi}{\rightarrow} & M \\
\iota'\,\uparrow\; & & \uparrow \iota \\
Z' & \overset{\psi'}{\rightarrow}& Z
\end{eqnarray*}
We define ${\psi'}^*$ with respect to these inclusions to smooth varieties. Then,
$\psi^*\iota_*[\mathcal F]=\iota'_*{\psi'}^*[\mathcal F]$. 
\end{quote}
To see this, observe that both sides are essentially the same 
$[\psi_*\mathcal O_N\otimes_{\mathcal O_M}^L\iota_*\mathcal F]$ 
by the definition of ${\psi'}^*$. 

Finally, recalling that $RR_n$ is defined by 
$$
RR_n(\mathcal F)=ch(\mathcal F)(1\otimes td_{Y^a_n})\cap
[Y_n^a\times Y_n^a],
$$
we have the commutativity in (2). 
\end{proof}

Theorem \ref{Borel-Moore} (1) also allows us to identify 
$H_{n-1,n}^a\rightarrow H_{n-1}^{a;i}$ with 
$$
H_*^{BM}(Z_{n-1,n}^a,\C)\longrightarrow H_*^{BM}(Z_{n-1}^{a;i},\C),
$$
but we do not try to find the geometric realization of this homomorphism. 
Later, $p_iH_{n-1,n}^ap_i\rightarrow H_{n-1}^{a;i}$ will be identified with 
$$
\Ext^*_{D^b(\Nil_n^a\times\mathbb P^m)}({R\pi_{n-1,n}^a}_!\C)
\longrightarrow \Ext^*_{D^b(\Nil_{n-1}^a)}({R\pi_{n-1}^{a;i}}_!\C),
$$
where the homomorphism is defined by the commutativity of the diagram: 
\begin{eqnarray*}
\Ext^*_{D^b(\Nil_n^a\times\mathbb P^m)}({R\pi_{n-1,n}^a}_!\C) &\longrightarrow
& \Ext^*_{D^b(\Nil_{n-1}^a)}({R\pi_{n-1}^{a;i}}_!\C)\\
\parallel\qquad\qquad  & & \quad\qquad\parallel \\
p_i\Hecke^a_{n-1,n}p_i\quad\quad & \longrightarrow & \;\quad\quad\Hecke^a_{n-1}
\end{eqnarray*}
Then, we consider the induced isomorphism 
$$
p_i\Hecke^a_{n-1,n}p_i/\Rad(p_i\Hecke^a_{n-1,n}p_i)
\simeq \Hecke^{a;i}_{n-1}/\Rad(\Hecke^{a;i}_{n-1})
$$
in this identification. Lemma \ref{functorial hom, part 2} shows that 
the isomorphism may be given geometrically and it suffices for our purpose. 

Recall that we have identified $e^\lambda\in R(T_n\times\C^\times)$ with 
$\prod_{i=1}^n X_i^{\lambda_i}\in \Hecke_n$. Denote the product by $X^\lambda$. 
Then, in the above theorem, $1\otimes X^\lambda$ is identified with  
$$
ch(\pi_n^*L_\lambda|_{Y^a_n})td_{Y^a_n}ch(\lambda_n)^{-1}\cap[Y^a_n]
\in H_*^{BM}(Y^a_n,\C).
$$
In particular, 
the identity element of $H_*^{BM}(Y^a_n,\C)$ is 
$td_{Y^a_n}ch(\lambda_n)^{-1}\cap[Y^a_n]$ 
and the multiplication by $X^\lambda$ is the same as the cap product 
$ch(\pi_n^*L_\lambda|_{Y^a_n})\cap-$.

\section{Geometric proof of the modular branching rule}

In this section, we give a geometric proof of the modular branching rule. 

\subsection{The statement}  
First we explain the precise statement which we are going to prove. In fact, 
we have two versions according to the choice of the identification 
$\Hecke_n=K^{G_n\times\C^\times}(Z_n)$. 

\begin{defn}
For an $\Hecke_n$-module $M$, define the $i$-restriction 
$$
i\text{-}\Res(M)=\{m\in M \mid (X_n-\zeta^i)^Nm=0, \text{for large enough $N$.}\}.
$$
\end{defn}

Then, the statement of the modular branching rule is as follows. The modules 
$L_\psi$ will be introduced in 4.4.

\begin{thm}
\label{The first branching theorem}
We identify $\Hecke_n$ with $K^{G_n\times\C^\times}(Z_n)$ by 
$\theta_\lambda=[{\delta_n}_*\pi_n^*L_\lambda]$ and 
$T_i=-[\mathcal L_i]-1$. Then, for the simple $\Hecke_n$-module $L_\psi$ 
labelled by an aperiodic multisegment $\psi$, we have
$$
\Soc(i\text{-}\Res(L_\psi))=L_{\tilde e_i\psi},
$$
where the crystal structure 
on the set of aperiodic multisegments is as in Theorem \ref{LTV-1}. 
\end{thm}

Let us consider the other identification of $\Hecke_n$ with 
$K^{G_n\times\C^\times}(Z_n)$. Recall the involution $\sigma$ 
defined by $T_i\mapsto -qT_i^{-1}$ and $X_i\mapsto X_i^{-1}$. 

\begin{defn}
An $\Hecke_n$-module obtained from $L_\psi$ by twisting the action by 
$\sigma$ and relabelling aperiodic multisegments by $\rho$ is denoted by
$$
D_\psi=\hphantom{}^\sigma L_{\rho(\psi)}.
$$
\end{defn}

\begin{thm}
\label{The second branching theorem}
We identify $\Hecke_n$ with $K^{G_n\times\C^\times}(Z_n)$ by 
$\theta_\lambda=[{\delta_n}_*\pi_n^*L_{-\lambda}]$ and 
$T_i=[\mathcal L_i]+q$. Then, for the simple $\Hecke_n$-module $D_\psi$ 
labelled by an aperiodic multisegment $\psi$, we have
$$
\Soc(i\text{-}\Res(D_\psi))=D_{\tilde e_i\psi}, 
$$
where the crystal structure 
on the set of aperiodic multisegments is as in Theorem \ref{LTV-2}. 
\end{thm}

Theorem \ref{The second branching theorem} follows from 
Theorem \ref{The first branching theorem}. In fact,  
\begin{equation*}
\begin{split}
\Soc(i\text{-}\Res(D_\psi))&\simeq 
\Soc(\hphantom{}^\sigma((-i)\text{-}\Res(L_{\rho(\psi)})))\\
&\simeq \hphantom{}^\sigma\Soc((-i)\text{-}\Res(L_{\rho(\psi)}))
\simeq \hphantom{}^\sigma L_{\tilde e_{-i}\rho(\psi)},
\end{split}
\end{equation*}
where $\tilde e_{-i}$ is the Kashiwara operator with respect to 
the crystal structure in Theorem \ref{The first branching theorem}, 
so that it is isomorphic to 
$\hphantom{}^\sigma L_{\rho(\tilde e_i\psi)}=D_{\tilde e_i\psi}$ where 
$\tilde e_i$ is the Kashiwara operator with respect to 
the crystal structure in Theorem \ref{The second branching theorem}. 

\textit{In the rest of the section, we identify $\Hecke_n$ with 
$K^{G_n\times\C^\times}(Z_n)$ as in Theorem 
\ref{The first branching theorem}} and prove the theorem. 

\subsection{Localization and eigenvalues of $X_n$}

Suppose that $(X,F)\in Y_n^a$. Then, $sXs^{-1}=\zeta X$ and $F$ is such that 
$F_i$ is obtained from $F_{i-1}$ by adding some  
eigenvector of $s$. We denote the eigenvalue of the eigenvector by 
$\zeta^{\nu_i}$, for $\nu_i\in\Z/e\Z$, and write $\nu=(\nu_1,\dots,\nu_n)$. 
We call $\nu$ the \emph{flag type} of $(X,F)$. Note that $\nu$ is a permutation of 
$(s_1,\dots,s_n)$. For $(X,F,F')\in Z_n^a=Y_n^a\times_{\Nil_n^a}Y_n^a$, we say that 
the \emph{flag type} of $(X,F,F')$ is $(\nu,\nu')$ if $(X,F)$ has flag type $\nu$ 
and $(X,F')$ has flag type $\nu'$. 

Now, we look at the decomposition of $Y_n^a$ and $Z^a_{n-1,n}$ 
into connected components. 
On each component, the flag type is constant. 

\begin{defn}
Let $p_iY_n^a$ be the disjoint union of connected components of 
$Y_n^a$ whose flag type $\nu$ satisfies $\nu_n=i$. 

Similarly, we let $p_iZ_{n-1,n}^ap_i$ be the disjoint union of 
connected components 
of $Z_{n-1,n}^a$ whose flag type $(\nu,\nu')$ satisfies $\nu_n=\nu'_n=i$. 
\end{defn}

The following lemma uses our choice of the identification of $\Hecke_n$ 
with $K^{G_n\times\C^\times}(Z_n)$ in this section. 

\begin{lem}
\label{eigenvalue analysis}
Under the identification 
$H_*^{BM}(Z_{n-1,n}^a,\C)=\Hecke_{n-1,n}^a$, we have 
$$
H_*^{BM}(p_iZ_{n-1,n}^ap_i,\C)=p_i\Hecke_{n-1,n}^ap_i. 
$$
\end{lem}
\begin{proof}
Let $(Y^a_n)_\mu$ be the set of 
$(X,F)\in Y^a_n$ such that the flag type is $\mu$. 
First we show that 
$$
H_*^{BM}(Y^a_n,\C)p_i=
\bigoplus_{\mu\;\text{such that}\;\mu_n=i}H_*^{BM}((Y^a_n)_\mu,\C).
$$
In fact, $X_n$ acts on $\C_a\otimes_{R(A)}K^A(Y^a_n)$ by
$$
\pi_n^*L_{\epsilon_n}|_{Y^a_n}\otimes -
$$
by Theorem \ref{Borel-Moore}. Now, $A$ acts on fiberwise over $Y^a_n$, and 
the fiber of $\pi_n^*L_{\epsilon_n}$ at $(X,F)$ is $\C^n/F_{n-1}$. 
Thus, $A$ acts as $\zeta^{\mu_n}$ on the fiber when the flag type of $(X,F)$ is $\mu$. 
Then, $X_n$ is 
$\zeta^{\mu_n}\pi_n^*L_{\epsilon_n}|_{Y^a_n}\in K((Y^a_n)_\mu)$,  
where $\pi_n^*L_{\epsilon_n}|_{Y^a_n}$ is a line bundle without $A$-action, 
and Theorem \ref{Borel-Moore} implies that 
$X_n$ acts on $H_*^{BM}((Y^a_n)_\mu,\C)$ by the cap product of 
$$
\zeta^{\mu_n}ch(\pi_n^*L_{\epsilon_n}|_{Y^a_n})=\zeta^{\mu_n}+\text{higher degree terms}.
$$
Hence, $X_n-\zeta^{\mu_n}$ acts nilpotently on $H_*^{BM}((Y^a_n)_\mu,\C)$. 
We have proved the claim. 

Let $\hphantom{}_{\nu}(Z^a_{n-1,n})_{\nu'}$ be the set of 
$(X,F,F')\in Z^a_{n-1,n}$ such that the flag type is $(\nu,\nu')$. 
By the definition of the convolution product, the product
$$
H_*^{BM}((Y^a_n)_\mu,\C)\cdot H_*^{BM}(\hphantom{}_{\nu}(Z^a_{n-1,n})_{\nu'},\C)
$$
is nonzero only if $\mu=\nu$. Thus, 
$p_iH_*^{BM}(\hphantom{}_{\nu}(Z^a_{n-1,n})_{\nu'},\C)=0$ if $\nu_n\ne i$, 
and the left multiplication by $p_i$ acts as the identity map on 
$H_*^{BM}(\hphantom{}_{\nu}(Z^a_{n-1,n})_{\nu'},\C)$ 
if $\nu_n=i$. Similar argument shows that 
$H_*^{BM}(\hphantom{}_{\nu}(Z^a_{n-1,n})_{\nu'},\C)p_i=0$ if $\nu'_n\ne i$, 
and the right multiplication by $p_i$ acts as the identity map on 
$H_*^{BM}(\hphantom{}_{\nu}(Z^a_{n-1,n})_{\nu'},\C)$ 
if $\nu'_n=i$. We have proved the result. 
\end{proof}

\subsection{A functorial algebra homomorphism} 
Now we work in the derived categories of abelian categories of 
sheaves of $\C$-vector spaces. 
The following is proved in \cite[Proposition 8.6.35]{CG}.

\begin{thm}
\label{sheaf theory}
Let $M_1$, $M_2$ and $M_3$ be connected smooth varieties, $N$ a variety and 
let $\mu_i: M_i\rightarrow N$ be proper maps. Let $\mathcal A_i\in D^b(M_i)$ be a constructible complex, for $i=1,2,3$. Define $Z_{ij}=M_i\times_N M_j$ and denote
$\iota_{ij}:Z_{ij}\subseteq M_i\times M_j$ the inclusion map. Let
$\mathcal A_{ij}=\iota_{ij}^!(\mathcal A_i^\vee\otimes \mathcal A_j)$. Then 
the following hold. 
\item[(1)]
Let $\mu_{ij}: Z_{ij}\rightarrow N$ be the projection map. Then
$$
R{\mu_{ij}}_*\mathcal A_{ij}\simeq 
R\mathcal Hom(R{\mu_i}_*\mathcal A_i,R{\mu_j}_*\mathcal A_j). 
$$ 
Thus, we have isomorphisms of $\C$-algebras 
$$
H^*(Z_{ij},\mathcal A_{ij})=H^*(N,R{\mu_{ij}}_*\mathcal A_{ij})\simeq 
\Ext^*_{D^b(N)}(R{\mu_i}_*\mathcal A_i,R{\mu_j}_*\mathcal A_j).
$$
\item[(2)]
The convolution product
$$
H^*(Z_{ij},\mathcal A_{ij})\otimes H^*(Z_{jk},\mathcal A_{jk})
\longrightarrow H^*(Z_{ik},\mathcal A_{ik})
$$
is identified with the Yoneda product
\begin{multline*}
\Ext^*_{D^b(N)}(R{\mu_i}_*\mathcal A_i,R{\mu_j}_*\mathcal A_j)\otimes 
\Ext^*_{D^b(N)}(R{\mu_j}_*\mathcal A_j,R{\mu_k}_*\mathcal A_k)\\
\longrightarrow 
\Ext^*_{D^b(N)}(R{\mu_i}_*\mathcal A_i,R{\mu_k}_*\mathcal A_k)
\end{multline*}
under the isomorphisms in (1). 
\end{thm}

We view elements of $\Nil^a_n$ as representations of the cyclic quiver of length 
$e$. Namely, we put $V_i=\{ v\in \C^n \mid sv=\zeta^iv \}$
on the $i^{th}$ node, for $i\in\Z/e\Z$, then $X\in\Nil^a_n$ defines 
$X: V_i\rightarrow V_{i+1}$, for $i\in\Z/e\Z$. 

We fix $i\in\Z/e\Z$. Let $m+1=\dim V_i$ and $\mathbb P^m$ the projective space 
consisting of $m$-dimensional subspaces of $V_i$. 
We have the following commutative diagram.
\begin{eqnarray*}
p_iY_n^a\underset{\Nil_n^a\times\mathbb P^m}{\times} p_iY_n^a=
p_iZ_{n-1,n}^ap_i & \hookrightarrow & \;\; 
p_iZ_n^ap_i=p_iY^a_n\underset{\Nil^a_n}{\times} p_iY^a_n \\
\downarrow\;\quad & & \quad\downarrow \\
\Nil_n^a\times \mathbb P^m & \underset{\rho_n}{\rightarrow} & \;\; \Nil_n^a 
\end{eqnarray*}
where $\rho_n(X,U)=X$ and the left vertical map is given by 
$(X,F,F')\mapsto(X,F_{n-1})$.

\begin{lem}
\label{functorial map}
Let $M \overset{f}{\rightarrow} X \overset{g}{\rightarrow} Y$ be proper maps and 
suppose that $M$ is smooth. We consider the following diagram, 
in which all squares are cartesian.
\begin{eqnarray*}
M\times_X M &\overset{\tilde\iota}{\longrightarrow} & M\times_Y M 
\overset{\tilde\Delta}{\longrightarrow} M\times M\\
\pi \downarrow & & \downarrow \pi' \qquad\qquad\;\; \downarrow f^{\times 2}\\
X & \overset{\iota}{\longrightarrow} & \;X\times_Y X 
\overset{\Delta}{\longrightarrow} X\times X\\
 & g\searrow\quad & \downarrow \pi'' \quad\qquad\swarrow g^{\times 2}\\
 &  & Y \overset{\overline\Delta}{\longrightarrow} Y \times Y
\end{eqnarray*}
Denote $\mathcal A=Rf_*\C$ and $\mathcal B=Rg_*\mathcal A$. 
Then the following hold.

\item[(1)]
We have the following isomorphisms of $\C$-algebras. 
$$
H_*^{BM}(M\times_X M,\C)\simeq \Ext^*_{D^b(X)}(\mathcal A,\mathcal A),\;\;
H_*^{BM}(M\times_Y M,\C)\simeq \Ext^*_{D^b(Y)}(\mathcal B,\mathcal B).
$$
\item[(2)]
$\tilde\iota_*: H_*^{BM}(M\times_X M,\C)\rightarrow H_*^{BM}(M\times_Y M,\C)$ 
is identified with the functorial algebra homomorphism
$$
Rg_*: \Ext^*_{D^b(X)}(\mathcal A,\mathcal A)\longrightarrow 
\Ext^*_{D^b(Y)}(\mathcal B,\mathcal B).
$$
\end{lem}
\begin{proof}
(1) follows from Theorem \ref{sheaf theory}. In fact, if we ignore degree shift then  
\begin{equation*}
\begin{split}
H_*^{BM}(M\times_X M,\C)&\simeq H^*(M\times_X M,\tilde\iota^!\tilde\Delta^!\C)
\simeq H^*(X, R\pi_*\tilde\iota^!\tilde\Delta^!\C)\\
&\simeq H^*(X, \iota^!R{\pi'}_*\tilde\Delta^!\C)
\simeq H^*(X, \iota^!\Delta^!Rf^{\times 2}_*\C).
\end{split}
\end{equation*}
As $\mathcal A^\vee=(Rf_*\C)^\vee=Rf_!\,\C^\vee=\oplus Rf_*\C[2\dim M_i]$, where 
the summation is over connected components $M_i$ of $M$, if we ignore degree shift then
$$
R\mathcal Hom_{D^b(X)}(\mathcal A,\mathcal A)=(\Delta\circ\iota)^!
(\mathcal A^\vee\otimes\mathcal A)=(\Delta\circ\iota)^!Rf^{\times 2}_*\C.
$$
Hence, $H_*^{BM}(M\times_X M,\C)\simeq \Ext^*_{D^b(X)}(\mathcal A,\mathcal A)$ is proved. 
We can prove the other isomorphism similarly.  

(2) If we ignore degree shift, the pushforward $\iota_*$ of 
Borel-Moore homology groups is given by
$$
H^*(M\times M, 
R(\tilde\Delta\circ\tilde\iota)_*(\tilde\Delta\circ\tilde\iota)^!(\C^\vee\otimes\C))
\longrightarrow
H^*(M\times M, R\tilde\Delta_*\tilde\Delta^!(\C^\vee\otimes\C)).
$$
First we claim that it is identified with 
$$
\Gamma(X\times_Y X, R\iota_*\iota^!\Delta^!(\mathcal A^\vee\otimes\mathcal A))
\longrightarrow
\Gamma(X\times_Y X, \Delta^!(\mathcal A^\vee\otimes\mathcal A)).
$$
To see this, let $\mathcal I^\bullet$ be an injective resolution of 
$\tilde\Delta^!(\C^\vee\otimes\C)$. Then, for the complex of sheaves 
$\Gamma_{M\times_X M}(\mathcal I^\bullet)$, which is defined by
$$
U\mapsto \Gamma_{M\times_X M}(\mathcal I^\bullet)(U)
=\{s^\bullet\in \mathcal I^\bullet(U) \mid \supp(s^i)\subseteq M\times_X M,\;
\text{for all $i$.} \},
$$
for open subsets $U\subseteq M\times_Y M$, the $\iota_*$ in question is obtained by 
taking the cohomology of the following morphism of complexes of $\C$-vector spaces. 
$$
\Gamma(M\times_Y M, \Gamma_{M\times_X M}(\mathcal I^\bullet))
\longrightarrow
\Gamma(M\times_Y M, \mathcal I^\bullet). 
$$
For open subsets $U\subseteq X\times_Y X$, we have
\begin{equation*}
\begin{split}
\Gamma_X(\pi'_*\mathcal F)(U)&=
\Ker\left(\pi'_*\mathcal F(U)\overset{restriction}{\longrightarrow}
\pi'_*\mathcal F(U\setminus X)\right)\\
&=
\Ker\left(\mathcal F({\pi'}^{-1}(U))\overset{restriction}{\longrightarrow}
\mathcal F({\pi'}^{-1}(U)\setminus M\times_X M)\right)\\
&=
\Gamma_{M\times_X M}(\mathcal F)({\pi'}^{-1}(U)),
\end{split}
\end{equation*}
for a sheaf $\mathcal F$ on $X\times_Y X$, so that the above morphism 
of complexes of $\C$-vector spaces is nothing but
$$
\Gamma(X\times_Y X, \Gamma_X(\pi'_*\mathcal I^\bullet))
\longrightarrow
\Gamma(X\times_Y X, \pi'_*\mathcal I^\bullet),
$$
and it is identified with  
$$
\Gamma(X\times_Y X, R\iota_*\iota^!R\pi'_*\tilde\Delta^!(\C^\vee\otimes\C))
\longrightarrow
\Gamma(X\times_Y X, R\pi'_*\tilde\Delta^!(\C^\vee\otimes\C)). 
$$
Now we apply the natural transformation $R\iota_*\iota^!\rightarrow \Id$ to 
the isomorphism
$$
R\pi'_*\tilde\Delta^!(\C^\vee\otimes\C)\simeq \Delta^!Rf_*^{\times 2}(\C^\vee\otimes\C)
\simeq \Delta^!((Rf_!\C)^\vee\otimes Rf_*\C))
$$
to obtain the claim. 

Next let $\mathcal I^\bullet$ be an injective resolution of $\mathcal A$. 
Then, our morphism of complexes of $\C$-vector spaces is 
$$
\Gamma(X\times_Y X, \Gamma_X({\mathcal I^\bullet}^\vee\otimes\mathcal I^\bullet))
\longrightarrow
\Gamma(X\times_Y X, 
\Gamma_{X\times_Y X}({\mathcal I^\bullet}^\vee\otimes\mathcal I^\bullet)). 
$$
For open subsets $U\subseteq X\times_Y X$, the map
$$
\Gamma_X({\mathcal I^\bullet}^\vee\otimes\mathcal I^\bullet)(U)
\longrightarrow
\Gamma_{X\times_Y X}({\mathcal I^\bullet}^\vee\otimes\mathcal I^\bullet)(U)
$$
sends $\sum \alpha_i^\bullet\otimes\beta_i^\bullet$, whose support is in $X$,  
to $\sum \alpha_i^\bullet\otimes\beta_i^\bullet$ itself. The left hand side 
is identified with 
$\Ext^*_{K^b(X)}(\mathcal I^\bullet, \mathcal I^\bullet)(U\cap X)$, where 
$K^b(X)$ is the homotopy category of the additive category of injective sheaves on $X$. 
On the other hand, if $U={\pi''}^{-1}(V)$, for an open subset $V\subseteq Y$, then 
$U\cap X=g^{-1}(V)$ and 
$$
\Gamma_{X\times_Y X}({\mathcal I^\bullet}^\vee\otimes\mathcal I^\bullet)(U)
=\Gamma_{\overline\Delta(Y)}(g_*{\mathcal I^\bullet}^\vee\otimes g_*\mathcal I^\bullet)(V)
$$
as before, so that the right hand side is identified with  
$\Ext^*_{K^b(Y)}(g_*\mathcal I^\bullet, g_*\mathcal I^\bullet)(V)$. 
Therefore, the pushforward $\iota_*$ of the Borel-Moore homology groups is 
the functorial algebra homomorphism $g_*$, namely $V=Y$ in the collection of maps
\begin{multline*}
g_*: \bigoplus_{i\in\Z}\Hom_{K^b(g^{-1}(V))}
(\mathcal I^\bullet|_{g^{-1}(V)}, \mathcal I^\bullet|_{g^{-1}(V)}[i])\\
\longrightarrow
\bigoplus_{i\in\Z}\Hom_{K^b(V)}(g_*\mathcal I^\bullet|_V, g_*\mathcal I^\bullet|_V[i]).
\end{multline*}
This is $Rg_*: \Ext^*_{D^b(X)}(\mathcal A, \mathcal A)\rightarrow 
\Ext^*_{D^b(Y)}(Rg_*\mathcal A, Rg_*\mathcal A)$ as desired. 
\end{proof}

In the following, we write $\Ext_{D^b(X)}^*(\mathcal A)$ for 
$\Ext_{D^b(X)}^*(\mathcal A,\mathcal A)$, and we denote 
\begin{align*}
\pi^a_{n-1,n}: p_iY^a_n &\longrightarrow \Nil^a_n\times\mathbb P^m,\\
\pi^a_n: p_iY^a_n &\longrightarrow \Nil^a_n.
\end{align*}
We remark that $\pi^a_n=\rho_n\circ\pi^a_{n-1,n}$. 

\begin{cor}
\label{functorial hom, part 1}
We have the isomorphisms 
\begin{align*}
p_i\Hecke_n^ap_i &\simeq\Ext^*_{D^b(\Nil_n^a)}({R\pi_n^a}_!\C),\\
p_i\Hecke_{n-1,n}^ap_i &\simeq 
\Ext^*_{D^b(\Nil_n^a\times\mathbb P^m)}({R\pi_{n-1,n}^a}_!\C)
\end{align*}
such that the inclusion 
$p_i\Hecke_{n-1,n}^ap_i\hookrightarrow p_i\Hecke_n^ap_i$ is identified with the 
following functorial algebra homomorphism. 
$$
R{\rho_n}_*: 
\Ext^*_{D^b(\Nil_n^a\times\mathbb P^m)}({R\pi_{n-1,n}^a}_!\C)
\longrightarrow\Ext^*_{D^b(\Nil_n^a)}({R\pi_n^a}_!\C).
$$
\end{cor}
\begin{proof}
Set $M=p_iY^a_n$, $X=\Nil^a_n\times\mathbb P^m$ and $Y=\Nil^a_n$. Then
Lemma \ref{functorial map} implies the result. 
\end{proof}

\subsection{Geometric construction of $U_v^-$} 
Let $U_v^-$ as in section 2. By Lusztig's theory, we may 
realize $U^-_v$ geometrically by using 
his geometric induction and restriction functors \cite{L2}. In fact, 
this is essentially the Hall algebra construction which we already 
explained in section 2. We only need 
the special case which corresponds to the multiplication by $f_i$, which 
we shall explain here. 

Recall that $\C^n$ has the eigenspace decomposition 
$\C^n=\oplus_{i\in\Z/e\Z}V_i$ with respect to 
$s=\diag(\zeta^{s_1},\dots,\zeta^{s_n})$. 
We suppose that $s_n=i$. 

Let $\dim V=m+1$ as before and 
let $W_i=V_i\cap\C^{n-1}$ and $W_j=V_j$, for $j\ne i$. 
Note that $W_i\ne V_i$. Thus we have 
$\C^{n-1}=\bigoplus_{i\in\Z/e\Z}W_i$. 
Then, we consider the diagram
$$
E_W\overset{p_1}{\longleftarrow} G_n(s)\times_{U_{n-1,n}(s)}F_{V,W}
\overset{p_2}{\longrightarrow} 
G_n(s)\times_{P_{n-1,n}(s)}F_{V,W}
\overset{p_3}{\longrightarrow} E_V,
$$
where $E_W$, $E_V$ and $F_{V,W}$ are defined by 
\begin{gather*}
E_W=\bigoplus_{i\in\Z/e\Z}\Hom_\C(W_i,W_{i+1}),\;\;
E_V=\bigoplus_{i\in\Z/e\Z}\Hom_\C(V_i,V_{i+1})\\
F_{V,W}=\{X\in E_V \mid XW_i\subseteq W_{i+1},\;\text{for all $i\in\Z/e\Z$.}\},
\end{gather*}
and $p_1(g,X)=X|_{\C^{n-1}}$, $p_2(g,X)=(g,X)$ and $p_3(g,X)=gXg^{-1}$. 

We only consider those objects whose supports are contained in the nullcones. 
We denote the subdiagram by
\footnote{Recall that $\Nil_{n-1,n}=\{ X\in \Nil_n \mid X\C^{n-1}\subseteq \C^{n-1} \}$.}
$$
\Nil^{a;i}_{n-1}\overset{\nu_{n-1,n}}{\longleftarrow}
G_n(s)\times_{U_{n-1,n}(s)}\Nil^a_{n-1,n}
\overset{\mu_{n-1,n}}{\longrightarrow}
G_n(s)\times_{P_{n-1,n}(s)}\Nil^a_{n-1,n}
\overset{\rho_n}{\longrightarrow}\Nil^a_n.
$$
Note that 
$G_n(s)\times_{P_{n-1,n}(s)}\Nil^a_{n-1,n}=
\{(X,U) \mid XU\subseteq U\}\subseteq \Nil^a_n\times\mathbb P^m$. 

$\Nil^{a;i}_{n-1}$ has finitely many $G_{n-1}(s)$-orbits and 
the stabilizer group of a point in each orbit $\mathcal O_\varphi$, 
for a multisegment $\varphi$, is connected. We denote by 
$$
IC_\varphi=IC(\overline{\mathcal O_\varphi},\C),
$$
the intersection cohomology complex associated with 
the orbit $\mathcal O_\varphi$ and the trivial local system on it. 
Then, $\nu_{n-1,n}^*IC_\varphi$ is a $L_{n-1,n}(s)$-equivariant simple perverse 
sheaf up to degree shift, and we may write 
$\nu_{n-1,n}^*IC_\varphi\simeq \mu_{n-1,n}^*IC_\varphi^\flat$ up to degree shift, 
for some simple perverse sheaf $IC_\varphi^\flat$ on $\Nil^a_n\times\mathbb P^m$. 
$IC_\varphi^\flat$ is unique up to isomorphism. In fact, 
we have an integer $d$ independent of $\varphi$, given by 
the difference of the dimensions 
of the fibers of $\mu_{n-1,n}$ and $\nu_{n-1,n}$, such that 
$IC_\varphi^\flat=
\hphantom{}^p\mathcal H^d({\nu_{n-1,n}}_*\mu_{n-1,n}^*IC_\varphi)$. 
We define a functor $\Ind_i^\flat$ by $\Ind_i^\flat(IC_\varphi)=IC_\varphi^\flat$. 
Then, we define the induction functor by
$$
\Ind_i={R\rho_n}_*\circ\Ind_i^\flat. 
$$
Now, as in the proof of \cite[9.2.3]{L2}, we consider the diagram
$$
Y^{a;i}_{n-1}\longleftarrow
G_n(s)\times_{U_{n-1,n}(s)}Y^a_{n-1,n}
\longrightarrow
G_n(s)\times_{P_{n-1,n}(s)}Y^a_{n-1,n}=p_iY^a_n,
$$
which \lq\lq covers\rq\rq\, the left three terms of the above diagram with cartesian squares. We denote the leftmost vertical map by
$$
\pi^{a;i}_{n-1}: Y^{a;i}_{n-1}\longrightarrow \Nil^{a;i}_{n-1}.
$$
Then, we have the following equalities up to degree shift. 
$$
\Ind_i^\flat({R\pi^{a;i}_{n-1}}_!\C)={R\pi^a_{n-1,n}}_!\C,\quad
\Ind_i({R\pi^{a;i}_{n-1}}_!\C)={R\pi^a_n}_!\C.
$$

The main result of \cite{L2} is the geometric construction of the 
algebra $U^-_v$ in terms of the induction functor.
The simple perverse sheaves $IC_\varphi$ 
are part of the canonical basis and $\Ind_i$ corresponds the multiplication 
from the left by $f_i$. 
The canonical basis defines the crystal $B(\infty)$. 
Combined with Kashiwara's result \cite[Proposition 6.2.3]{Ka}, 
we have the following. \footnote{It is known that \cite[Proposition 6.2.3]{Ka} 
may be proved in this geometric framework.}

\begin{lem}
\label{degree control}
\item[(1)]
Let $\varphi$ be a multisegment of size $n-1$. Then, we may write
$$
\Ind_i(IC_\varphi)=
\bigoplus_{j=0}^{\epsilon_i(\varphi)}
IC_{\tilde f_i\varphi}[\epsilon_i(\varphi)-2j]\bigoplus
\left(\bigoplus_{j\in\Z} R_{\varphi,j}[j]\right),
$$
for certain perverse sheaves $R_{\varphi,j}$ on $\Nil^a_n$. 
\item[(2)]
Suppose that $IC_\psi$, for a multisegment $\psi$ of size $n$, 
appears in $R_{\varphi,j}$, for some $j$. Then, we have 
$$
-\epsilon_i(\psi)+2\le j\le \epsilon_i(\psi)-2.
$$
\end{lem}

\subsection{Some semisimple quotients} 
Lemma \ref{simple modules} implies that the surjection 
$p_i\Hecke^a_{n-1,n}p_i\rightarrow\Hecke^{a;i}_{n-1}$ 
induces the identity map
\begin{multline*}
p_i\Hecke^a_{n-1,n}p_i/\Rad(p_i\Hecke^a_{n-1,n}p_i)\simeq\oplus_M \Endhom_\C(M)\\
\longrightarrow \oplus_M \Endhom_\C(M)\simeq\Hecke^{a;i}_{n-1}/\Rad(\Hecke^{a;i}_{n-1}), 
\end{multline*}
where $M$ runs through the common complete set of isomorphism classes of simple modules. 

On the other hand, 
the complete set of isomorphism classes of simple modules of 
$\Ext^*_{D^b(\Nil_n^a\times\mathbb P^m)}({R\pi_{n-1,n}^a}_!\C)$ and 
$\Ext^*_{D^b(\Nil_{n-1}^{a;i})}({R\pi_{n-1}^{a;i}}_!\C)$ may be described by 
those simple perverse sheaves that appear in 
${R\pi_{n-1,n}^a}_!\C$ and ${R\pi_{n-1}^{a;i}}_!\C$ after some shift, respectively. 
The degree of the shift depends on the perverse sheaf.  
As they are semisimple complexes by the decomposition theorem, we write
$$
{R\pi_{n-1,n}^a}_!\C\simeq\sum_{\psi}\sum_{m\in\Z} IC_\psi[m]^{\oplus m_{\psi,m}},\quad
{R\pi_{n-1}^{a;i}}_!\C\simeq\sum_{\varphi}\sum_{m\in\Z} IC_\varphi[m]^{\oplus n_{\varphi,m}},
$$
where $IC_\psi$ and $IC_\varphi$ are simple perverse sheaves on 
$\Nil_n^a\times\mathbb P^m$ and $\Nil_{n-1}^{a;i}$, respectively. 
Let $L_{\psi,m}=\C^{m_{\psi,m}}$ and $L_{\varphi,m}=\C^{m_{\varphi,m}}$ 
be the multiplicity spaces of $IC_\psi[m]$ and $IC_\varphi[m]$, respectively.
Define 
$$
L_\psi=\bigoplus_{m\in\Z}L_{\psi,m},\quad L_\varphi=\bigoplus_{m\in\Z}L_{\varphi,m}.
$$
Then, we have
\begin{equation*}
\begin{split}
\Ext^*_{D^b(\Nil_n^a\times\mathbb P^m)}({R\pi_{n-1,n}^a}_!\C)&\simeq
\bigoplus_{\psi',\psi''}
\Ext^*(IC_{\psi'},IC_{\psi''})\otimes_\C\Hom_\C(L_{\psi'},L_{\psi''}),\\
\Ext^*_{D^b(\Nil_{n-1}^{a;i})}({R\pi_{n-1}^{a;i}}_!\C)&\simeq
\bigoplus_{\varphi',\varphi''}
\Ext^*(IC_{\varphi'},IC_{\varphi''})\otimes_\C\Hom_\C(L_{\varphi'},L_{\varphi''}).
\end{split}
\end{equation*}
In other words, $\Ext^*_{D^b(\Nil_n^a\times\mathbb P^m)}({R\pi_{n-1,n}^a}_!\C)$ 
is the matrix algebra which has block partitions of rows and columns such that 
the blocks are labelled by $\psi$ and the entries in the $(\psi'',\psi')$ component
are elements of $\Ext^*(IC_{\psi'},IC_{\psi''})$. In particular, its 
largest semisimple quotient is the 
block diagonal matrix algebra such that the entries of the $(\psi,\psi)$-component 
are
$$
\Ext^{\ge0}(IC_\psi,IC_\psi)/\Ext^{>0}(IC_\psi,IC_\psi)\simeq\C.
$$ 
We have the similar matrix 
algebra description for $\Ext^*_{D^b(\Nil_{n-1}^{a;i})}({R\pi_{n-1}^{a;i}}_!\C)$ as well. 

\subsection{A key result} 
We prove Theorem \ref{key result}, which we will need in the geometric proof 
of the modular branching rule in the next subsection. 

Define $\eta_{n-1,n}: \Nil^a_{n-1,n}\hookrightarrow \Nil^a_n\times\mathbb P^m$,  
$\kappa_{n-1,n}: \Nil^a_{n-1,n}\rightarrow \Nil^{a;i}_{n-1}$. We identify 
$\Nil^{a;i}_{n-1}$ with the zero section of $\kappa_{n-1,n}$ and we obtain 
the closed embedding
$$
\epsilon_{n-1,n}: \Nil^{a;i}_{n-1}\hookrightarrow \Nil^a_n\times\mathbb P^m.
$$ 
$\eta_{n-1,n}^*{R\pi_{n-1,n}^a}_!\C$ is the pushforward of the constant sheaf 
on $Y^a_{n-1,n}$ to $\Nil^a_{n-1,n}$, and we have the following cartesian diagram.
\begin{eqnarray*}
Y^a_{n-1,n} & \rightarrow & Y^{a;i}_{n-1}\\
\downarrow\quad & & \quad\downarrow \\
\Nil^a_{n-1,n} & \rightarrow & \Nil^{a;i}_{n-1}
\end{eqnarray*}
Thus, $\eta_{n-1,n}^*{R\pi_{n-1,n}^a}_!\C\simeq 
\kappa_{n-1,n}^*{R\pi_{n-1}^{a;i}}_!\C$ and we conclude that 
$$
\epsilon_{n-1,n}^*{R\pi_{n-1,n}^a}_!\C\simeq {R\pi_{n-1}^{a;i}}_!\C. 
$$

\begin{lem}
\label{functorial hom, part 2}
We consider the functorial algebra homomorphism
$$
\epsilon_{n-1,n}^*: \Ext^*_{D^b(\Nil_n^a\times\mathbb P^m)}({R\pi_{n-1,n}^a}_!\C)
\longrightarrow \Ext^*_{D^b(\Nil_{n-1}^{a;i})}({R\pi_{n-1}^{a;i}}_!\C).
$$
Then, it induces the isomorphism
\begin{gather*}
\Ext^*_{D^b(\Nil_n^a\times\mathbb P^m)}({R\pi_{n-1,n}^a}_!\C)/
\Rad(\Ext^*_{D^b(\Nil_n^a\times\mathbb P^m)}({R\pi_{n-1,n}^a}_!\C))\\
\qquad\simeq
\Ext^*_{D^b(\Nil_{n-1}^{a;i})}({R\pi_{n-1}^{a;i}}_!\C)/
\Rad(\Ext^*_{D^b(\Nil_{n-1}^{a;i})}({R\pi_{n-1}^{a;i}}_!\C))
\end{gather*}
and it is identified with the identity map 
$$
p_i\Hecke^a_{n-1,n}p_i/\Rad(p_i\Hecke^a_{n-1,n}p_i)
\simeq \Hecke^{a;i}_{n-1}/\Rad(\Hecke^{a;i}_{n-1}).
$$
Further, its inverse is induced by the functorial algebra homomorphism
$$
\Ind_i^\flat: \Ext^*_{D^b(\Nil_{n-1}^{a;i})}({R\pi_{n-1}^{a;i}}_!\C)
\longrightarrow \Ext^*_{D^b(\Nil_n^a\times\mathbb P^m)}({R\pi_{n-1,n}^a}_!\C).
$$
\end{lem}
\begin{proof}
Note that 
$$
\Ext^*_{D^b(\Nil_{n-1}^{a;i})}({R\pi_{n-1}^{a;i}}_!\C)=
\Ext^*_{D^b(\Nil^{a;i}_{n-1})}(\oplus_\varphi IC_\varphi\otimes_\C L_\varphi)
$$
as $\C$-algebras. Thus, the functorial algebra homomorphism
$$
\Ind_i^\flat: \Ext^*_{D^b(\Nil^{a;i}_{n-1})}(\oplus_\varphi IC_\varphi\otimes_\C L_\varphi)
\rightarrow
\Ext^*_{D^b(\Nil^a_{n-1,n})}(\oplus_\varphi IC_\varphi^\flat\otimes_\C L_\varphi)
$$
induces the identity map 
\begin{multline*}
\Ext^0_{D^b(\Nil^{a;i}_{n-1})}(\oplus_\varphi IC_\varphi\otimes_\C L_\varphi)=
\oplus_\varphi \Endhom_\C(L_\varphi)\\
\longrightarrow
\oplus_\varphi \Endhom_\C(L_\varphi)=
\Ext^0_{D^b(\Nil^a_{n-1,n})}(\oplus_\varphi IC_\varphi^\flat\otimes_\C L_\varphi).
\end{multline*}
That is, $\Ind_i^\flat$ induces the isomorphism
\begin{gather*}
\Ext^*_{D^b(\Nil_n^a\times\mathbb P^m)}({R\pi_{n-1}^{a;i}}_!\C)/
\Rad(\Ext^*_{D^b(\Nil_n^a\times\mathbb P^m)}({R\pi_{n-1}^{a;i}}_!\C))\\
\qquad\simeq
\Ext^*_{D^b(\Nil_{n-1}^{a;i})}({R\pi_{n-1,n}^a}_!\C)/
\Rad(\Ext^*_{D^b(\Nil_{n-1}^{a;i})}({R\pi_{n-1,n}^a}_!\C)).
\end{gather*}
it is identified with the identity map 
$$
\Hecke^{a;i}_{n-1}/\Rad(\Hecke^{a;i}_{n-1})\simeq 
p_i\Hecke^a_{n-1,n}p_i/\Rad(p_i\Hecke^a_{n-1,n}p_i). 
$$
On the other hand, we have 
$\Ind_i^\flat({R\pi_{n-1}^{a;i}}_!\C)\simeq {R\pi_{n-1,n}^a}_!\C$ 
up to degree shift, and 
$\epsilon_{n-1,n}^*{R\pi_{n-1,n}^a}_!\C\simeq {R\pi_{n-1}^{a;i}}_!\C$. 
Thus, $\Ind_i^\flat$ and $\epsilon_{n-1,n}^*$ are inverse to the other on 
the semisimple quotients, and the claim follows. 
\end{proof}

\begin{thm}
\label{key result}
Consider the functorial algebra homomorphism
$$
\Ind_i: \Ext_{D^b(\Nil_{n-1}^{a;i})}^*({R\pi_{n-1}^{a;i}}_!\C)\longrightarrow
\Ext_{D^b(\Nil_n^a)}^*({R\pi_n^a}_!\C).
$$
If $M$ is a simple $\Hecke^a_n$-module, then 
the action of $\Hecke^{a;i}_{n-1}$ on $\Top(p_iM)$ coincides with that 
given by $\Ind_i$ under the identification
$$
\Hecke^{a;i}_{n-1}=\Ext_{D^b(\Nil_{n-1}^{a;i})}^*({R\pi_{n-1}^{a;i}}_!\C),\;\;
p_i\Hecke^a_np_i=\Ext_{D^b(\Nil_n^a)}^*({R\pi_n^a}_!\C).
$$
\end{thm}
\begin{proof}
Let $(Y_n^a)_\nu$ be the set of 
$(X,F)$ such that the flag type is $\nu$, as before. We denote
$\pi_{n,\nu}: (Y_n^a)_\nu\rightarrow \Nil^a_n$ and 
$$
\mathcal M_\nu=\bigoplus_{i\in\Z}\hphantom{}^p\mathcal H^i(R{\pi_{n,\nu}}_!\C). 
$$
Then, by our identification, we have
$$
\Hecke_n^a=\Ext^*_{D^b(\Nil_n^a)}(\oplus_\nu \mathcal M_\nu)
$$
where $\nu$ runs through flag types which are permutations of 
$(s_1,\dots,s_n)$. Write
$$
\bigoplus_\nu \mathcal M_\nu=\bigoplus_\psi IC_\psi\otimes_\C L_\psi.
$$
Then,
$\Hecke^a_n=\bigoplus_{\psi',\psi''} 
\Ext^*_{D^b(\Nil_n^a)}(IC_{\psi'}, IC_{\psi''})
\otimes_\C \Hom_\C(L_{\psi'}, L_{\psi''})$ 
and we view it as the block partitioned matrix algebra whose 
entries of the $(\psi'',\psi')$-component are elements of 
$\Ext^*_{D^b(\Nil_n^a)}(IC_{\psi'}, IC_{\psi''})$. Define
$$
P_\psi=\bigoplus_{\psi'} \Ext^*_{D^b(\Nil_n^a)}(IC_\psi, IC_{\psi'})
\otimes_\C L_{\psi'}.
$$
Then, it is a direct summand of $\Hecke_n^a$ and we view it 
as the space of block partitioned
column vectors whose entries in the block $L_{\psi'}$ are elements of
$\Ext^*_{D^b(\Nil_n^a)}(IC_\psi, IC_{\psi'})$. 
\begin{multline*}
\Ext^*_{D^b(\Nil_n^a)}(IC_{\psi'}, IC_{\psi''})\otimes_\C\Hom_\C(L_{\psi'}, L_{\psi''})
\times \Ext^*_{D^b(\Nil_n^a)}(IC_\psi, IC_{\psi'})\otimes_\C L_{\psi'}\\
\longrightarrow
\Ext^*_{D^b(\Nil_n^a)}(IC_\psi, IC_{\psi''})\otimes_\C L_{\psi''}
\end{multline*}
shows that $P_\psi$ is a left ideal of $\Hecke_n^a$ 
so that it is a projective $\Hecke_n^a$-module. It is clear that 
$$
L_\psi=\frac{\Ext^{\ge0}_{D^b(\Nil_n^a)}(IC_\psi, \oplus_\nu \mathcal M_\nu)}
{\Ext^{>0}_{D^b(\Nil_n^a)}(IC_\psi, \oplus_\nu \mathcal M_\nu)}
$$
is a simple $\Hecke_n^a$-module or zero and that any simple $\Hecke_n^a$-module appears 
in this way. Thus, we assume that $M=L_\psi$. 
Then, Lemma \ref{eigenvalue analysis} says that multiplication by $p_i$ amounts 
to picking out the connected components $p_iY_n^a$ so that
$$
p_iL_\psi=\frac{\Ext^{\ge0}_{D^b(\Nil_n^a)}(IC_\psi, \oplus_\nu \mathcal M_\nu)}
{\Ext^{>0}_{D^b(\Nil_n^a)}(IC_\psi, \oplus_\nu \mathcal M_\nu)}
$$
where $\nu$ runs through permutations of $(s_1,\dots,s_n)$ such that 
$\nu_n=i$. Suppose that $p_iL_\psi\ne0$. It is a simple 
$p_i\Hecke^a_np_i$-module. Let 
$\pi_{n-1,n,\nu}: (Y_n^a)_\nu\rightarrow \Nil^a_n\times\mathbb P^m$ and 
$$
\mathcal M_\nu^\flat=\bigoplus_{i\in\Z}
\hphantom{}^p\mathcal H^i(R{\pi_{n-1,n,\nu}}_!\C). 
$$
Then $p_i\Hecke^a_{n-1,n}p_i=
\Ext_{D^b(\Nil^a_n\times\mathbb P^m)}(\oplus_\nu \mathcal M_\nu^\flat)$, 
where $\nu$ runs through permutations of $(s_1,\dots,s_n)$ such that 
$\nu_n=i$, and it acts on $p_iL_\psi$ through ${R\rho_n}_*$ by 
Corollary \ref{functorial hom, part 1}. Now, we consider 
$\Top(p_iL_\psi)$. Then, the action of $p_i\Hecke^a_{n-1,n}p_i$ factors 
through $\Hecke^{a;i}_{n-1}/\Rad(\Hecke^{a;i}_{n-1})$ and 
Lemma \ref{functorial hom, part 2} implies that it is given by 
$\Ind_i^\flat$. We have proved that the action of $\Hecke^{a;i}_{n-1}$ on 
$\Top(p_iL_\psi)$ coincides with the action of $\Hecke^{a;i}_{n-1}$ given by 
the functorial algebra homomorphism $\Ind_i$. 
\end{proof}

\subsection{The geometric proof} 

Having proved Theorem \ref{key result}, we are now able to give the 
promised geometric proof 
of the modular branching rule. We write each simple $\Hecke^a_n$-module as 
in the proof of the above theorem
$$
L_\psi=\frac{\Ext^{\ge0}_{D^b(\Nil_n^a)}(IC_\psi, \oplus_\nu \mathcal M_\nu)}
{\Ext^{>0}_{D^b(\Nil_n^a)}(IC_\psi, \oplus_\nu \mathcal M_\nu)}.
$$
Suppose that $p_iL_\psi\ne0$. We want to show that $\Top(p_iL_\psi)$ 
contains $L_{\tilde e_i\psi}$. As the simple $\Hecke^{a;i}_{n-1}$-modules 
are the same as the simple $p_i\Hecke^a_{n-1,n}p_i$-modules, we consider 
the restriction of $p_iL_\psi$ to $p_i\Hecke^a_{n-1,n}p_i$. Let 
$\pi^a_{n,\nu}=\rho_n\circ\pi^a_{n-1,n,\nu}$. Then, we have
$$
{R\pi^a_n}_!\C=\bigoplus_{\text{$\nu$ such that $\nu_n=i$}}
{R\pi^a_{n,\nu}}_!\C,
$$
which is equal to $\Ind_i({R\pi^{a;i}_{n-1}}_!\C)$ up to degree shift. 
Thus, we write
$$
\bigoplus_{\text{$\nu$ such that $\nu_n=i$}}
\mathcal M_\nu^\flat=\bigoplus_{\varphi}
IC_\varphi^\flat\otimes_\C L_\varphi
$$
and restrict the action of $p_i\Hecke^a_np_i$ on 
$p_iL_\psi$ to $p_i\Hecke^a_{n-1,n}p_i$ through $R{\rho_n}_*$, which is 
the functorial algebra homomorphism given by 
\begin{multline*}
p_i\Hecke^a_{n-1,n}p_i=\bigoplus_{\varphi',\varphi''}
\Ext^*_{D^b(\Nil_n^a\times\mathbb P^m)}
(\Ind_i^\flat IC_{\varphi'}, \Ind_i^\flat IC_{\varphi''})
\bigotimes_{\C}\Hom_\C(L_{\varphi'}, L_{\varphi''})\\
\longrightarrow
\bigoplus_{\varphi',\varphi''}
\Ext^*_{D^b(\Nil_n^a)}(\Ind_i IC_{\varphi'}, \Ind_i IC_{\varphi''})
\bigotimes_{\C}\Hom_\C(L_{\varphi'}, L_{\varphi''})=p_i\Hecke^a_np_i.
\end{multline*}

To study this, we introduce a block algebra description of 
$p_i\Hecke^a_{n-1,n}p_i$-action 
on $p_iL_\psi$. As
$$
\bigoplus_{\text{$\nu$ such that $\nu_n=i$}}\mathcal M_\nu
=\bigoplus_\varphi
\left(IC_{\tilde f_i\varphi}^{\oplus(\epsilon_i(\varphi)+1)}
+\sum_j R_{\varphi,j}\right)\otimes_\C L_\varphi,
$$
by Lemma \ref{degree control}(1), $p_iL_\psi$ has the decomposition
$$
p_iL_\psi=
\bigoplus_\varphi\frac{\Ext^{\ge0}_{D^b(\Nil_n^a)}
(IC_\psi, IC_{\tilde f_i\varphi}^{\oplus(\epsilon_i(\varphi)+1)}
+\sum_j R_{\varphi,j})}
{\Ext^{>0}_{D^b(\Nil_n^a)}
(IC_\psi, IC_{\tilde f_i\varphi}^{\oplus(\epsilon_i(\varphi)+1)}
+\sum_j R_{\varphi,j})}\otimes_\C L_\varphi.
$$
Thus, we have the corresponding block decomposition of $\Endhom_\C(p_iL_\psi)$. 

Observe that $IC_\psi$ appears in $R_{\varphi,j}$ only if 
$\epsilon_i(\varphi)<\epsilon_i(\tilde e_i\psi)$ and 
$IC_\psi$ appears in $IC_{\tilde f_i\varphi}^{\oplus(\epsilon_i(\varphi)+1)}$ 
only if $\varphi=\tilde e_i\psi$. Hence, only those blocks 
$L_\varphi$ with $\epsilon_i(\varphi)<\epsilon_i(\tilde e_i\psi)$ and 
$L_{\tilde e_i\psi}$ appear in the above block decomposition. 

To obtain the $(\varphi'', \varphi')$-component of 
the representation of $p_i\Hecke^a_{n-1,n}p_i$ on $p_iL_\psi$, 
we consider the image of 
$\Ext^k_{D^b(\Nil^a_n\times\mathbb P^m)}(IC_{\varphi'}^\flat, IC_{\varphi''}^\flat)$, 
for $k\ge0$, through the action of 
$$
\Ext^k_{D^b(\Nil^a_n)}(\Ind_i(IC_{\varphi'}), \Ind_i(IC_{\varphi''})).
$$ 
The image may be nonzero only when $IC_\psi[j']$, for some $j'\in\Z$, appears 
in $\Ind_i(IC_{\varphi'})$ and $IC_\psi[j'']$, for some $j''\in\Z$, appears 
in $\Ind_i(IC_{\varphi''})$ such that $-j'+j''+k=0$. In particular, 
$j''\le j'$ is necessary. Since
\begin{equation*}
\begin{split}
\Ind_i(IC_{\varphi'})&=
\sum_{j'=0}^{\epsilon_i(\varphi')}IC_{\tilde f_i\varphi'}[\epsilon_i(\varphi')-2j']
+\sum_{j'\in\Z}R_{\varphi',j'}[j'],\\
\Ind_i(IC_{\varphi''})&=
\sum_{j''=0}^{\epsilon_i(\varphi'')}IC_{\tilde f_i\varphi''}[\epsilon_i(\varphi'')-2j'']
+\sum_{j''\in\Z}R_{\varphi'',j''}[j''],
\end{split}
\end{equation*}
there are four cases to consider. 
\begin{itemize}
\item
Suppose that $\varphi'=\varphi''=\tilde e_i\psi$. We number the rows and columns 
of the block matrix by $0\le j'', j'\le\epsilon_i(\psi)-1$ such that 
$\epsilon_i(\psi)-1-2j''$ and $\epsilon_i(\psi)-1-2j'$ are increasing. 
Then, the entries may be 
nonzero only when $\epsilon_i(\psi)-1-2j''\le \epsilon_i(\psi)-1-2j'$. Thus, 
we obtain an upper block triangular matrix whose diagonal block components are 
$\Endhom_\C(L_{\tilde e_i\psi})$. 
\item
Suppose that $\varphi'\ne\tilde e_i\psi=\varphi''$. We number the rows as before, 
and the columns such that $j'$ is increasing. If $IC_\psi$ appears in 
$L_{\varphi',j'}$ then the entries may be 
nonzero only when $\epsilon_i(\psi)-1-2j''\le j'$. Hence, each row 
has entries only after the column number $\epsilon_i(\psi)-1-2j''$. 
Now, Lemma \ref{degree control}(2) implies that $j'\le \epsilon_i(\psi)-2$ 
so that $j''=0$ cannot happen. Hence, all the entries of the last row are zero. 
\item
Suppose that $\varphi'=\tilde e_i\psi\ne\varphi''$. Then, 
each column has entries only before some column number. 
\item
Suppose that $\varphi'\ne\tilde e_i\psi$ and $\varphi''\ne\tilde e_i\psi$. 
Then we have an upper block triangular matrix again. 
\end{itemize}
The first two cases show that there is a $p_i\Hecke^a_{n-1,n}p_i$-submodule 
$L'_\psi$ of $L_\psi$ such that $L_\psi/L'_\psi\simeq L_{\tilde e_i\psi}$. 
Thus, $L_{\tilde e_i\psi}$ appears in $\Top(p_iL_\psi)$. 
Now, following \cite{Kl}, Grojnowski and Vazirani proved in Vazirani's thesis that 
$\Soc(p_iL_\psi)$ is simple \cite{GV}. By Specht module theory, 
the simple modules are self-dual so that  $\Top(p_iL_\psi)$ is 
isomorphic to $\Soc(p_iL_\psi)$. Thus, we have proved that 
$\Soc(p_iL_\psi)=L_{\tilde e_i\psi}$. Thus, 
Theorem \ref{The first branching theorem} and 
Theorem \ref{The second branching theorem} follow. 

\section{Crystals of deformed Fock spaces}

In this section, we recall results on deformed Fock spaces which are related to 
the combinatorial construction of simple $\Hecke_n$-modules. 

\subsection{Crystals of deformed Fock spaces} 
Let $l\in \Z_{>0}$ and we choose a \emph{multicharge} 
$$
\mathbf v=(\mathrm v_0,...,\mathrm v_{l-1})\in \Z^l. 
$$
We denote $\mathrm v_i+e\Z\in \Z/e\Z$ by $\overline{\mathrm v_i}$, for $1\le i\le l$. 
Let $\Lambda_i$, for $i\in\Z/e\Z$, be the fundamental weights of $\mathfrak g$, 
and define a dominant weight $\Lambda$ by 
$$
\Lambda=\Lambda_{\overline{\mathrm v_0}}+\cdots+\Lambda_{\overline{\mathrm v_{l-1}}}. 
$$
We consider various multicharges which give a fixed $\Lambda$. 

Let $V_v(\Lambda)$ be the integrable highest weight $U_v(\mathfrak g)$-module 
of highest weight $\Lambda$. We want to realize $V_v(\Lambda)$ as a 
$U_v(\mathfrak g)$-submodule of 
the level $l$ deformed Fock space $\Fock^{\mathbf v}$ 
associated with the multicharge $\mathbf v$. 

As a $\C(v)$-vector space, the level $l$ Fock space $\Fock^{\mathbf v}$
admits the set of all $l$-partitions as a natural basis. Namely, 
the underlying vector space is 
$$
\Fock=\bigoplus_{n\geq 0}\bigoplus_{\ulambda\in \Pi_{l,n}}
\C(v)\ulambda,
$$
where $\Pi_{l,n}$ is the set of $l$-partitions of rank $n$. We do not 
give explicit formulas to define 
the $U_v(\mathfrak g)$-module structure on $\Fock^{\mathbf v}$, 
but it is defined in terms of the total order $\prec_{\mathbf v}$ 
introduced below. 
The action we adopt here is the one which was introduced by Jimbo, Misra, Miwa and Okado in \cite{jim}. 
Let 
$$
L^{\mathbf v}=\bigoplus_{n\geq 0}\bigoplus_{\ulambda\in \Pi_{l,n}}R\ulambda,\quad 
B^{\mathbf v}=\bigsqcup_{n\geq 0} \Pi_{l,n}.
$$
Then, $(L^{\mathbf v}, B^{\mathbf v})$ is a crystal basis of $\Fock^{\mathbf v}$. 
In this article, it suffices to recall the crystal structure on the 
set of $l$-partitions. Before doing this, we explain basic terminology 
on $l$-partitions.  

Let $\ulambda=(\lambda^{(0)},\dots,\lambda^{(l-1)})$ be an $l$-partition, 
which is identified with the corresponding $l$-tuple of Young diagrams. 
Then, we can speak of nodes of $\ulambda$, which are nodes of 
the Young diagrams. We identify a node $\gamma$ of $\ulambda$ with a triplet  
$(a,b,c)$ where $c\in \{0,...,l-1\}$ is such that $\gamma$ is a node of 
$\lambda^{(c)}$, and $a$ and $b$ are the
row and the column indices of the node $\gamma$ in $\lambda^{(c)}$, respectively. 

\begin{defn}
Let $\gamma=(a,b,c)$ be a node of an $l$-partition $\ulambda$. Then, 
its \emph{content} $c(\gamma)$ and \emph{residue} 
$\residue(\gamma)$ are defined by
$$
c(\gamma)=b-a+\mathrm v_c\in\Z\;\;\text{and}\;\;
\residue(\gamma)=\overline{c(\gamma)}\in\Z/e\Z,
$$
respectively. 
\end{defn}

Let $\gamma$ be a node of $\ulambda$. Then we say that 
$\gamma$ is an \emph{$i$-node}, for $i\in\Z/e\Z$, if 
$\residue(\gamma)=i$. Suppose that 
$\ulambda\backslash \{\gamma\}$ is again an $l$-partition, 
which we denote by $\umu$. Then, we say that 
$\gamma$ is a \emph{removable $i$-node} 
of $\ulambda$ and $\gamma$ is an \emph{addable $i$-node} 
of $\umu$. We introduce a total order $\prec_{\mathbf v}$ 
on the set of addable and removable $i$-nodes
of an $l$-partition $\ulambda$. 

\begin{defn}
Let $\gamma_1=(a_1,b_1,c_1)$ and $\gamma_2=(a_2,b_2,c_2)$ be 
$i$-nodes of $\ulambda$.
We define the order $\prec_{\mathbf v}$ by 
$$
\gamma _{1}\prec _{{\mathbf{v}}}\gamma _{2}\Longleftrightarrow 
\begin{cases}
c(\gamma _{1})<c(\gamma _{2}),\text{ or} \\ 
c(\gamma _{1})=c(\gamma _{2})\text{ and }c_{1}>c_{2}.
\end{cases}
$$
\end{defn}
The order $\prec_{\mathbf v}$ depends on the choice of 
the multicharge $\mathbf v$ when $l>1$. 

Now, we can explain the crystal structure on $B^{\mathbf v}$, 
which is defined by the total order $\prec_{\mathbf v}$. 
Let $\ulambda$ be an $l$-partition as above. 
Let $N_i(\ulambda)$, for $i\in\Z/e\Z$, be 
the number of $i$-nodes of $\ulambda$. Then we define
$$
\wt(\ulambda)=\Lambda-\sum_{i\in\Z/e\Z}N_i(\ulambda)\alpha_i.
$$
The rule 
to compute $\tilde e_i\ulambda$ is as follows. 
The rule to compute $\tilde f_i\ulambda$ is similar. 

\begin{quote}
We read addable and removable $i$-nodes of $\ulambda$ 
in the increasing order with respect to $\prec_{\mathbf v}$. 
Then we delete a consecutive pair of a removable $i$-node and 
an addable $i$-node in this order as many as possible. 
We call this procedure $RA$ deletion. 
\begin{itemize}
\item
If there remains no removable $i$-node, 
define $\tilde e_i\ulambda=0$. 
\item
Otherwise, we call the leftmost removable $i$-node, say $\gamma$, 
the \emph{good $i$-node of $\ulambda$}, and define 
$\tilde e_i\ulambda=\ulambda\backslash\{\gamma\}$.
\end{itemize} 
\end{quote}
Finally, we define 
$$
\epsilon_i(\ulambda)=\max\{k\in\Z_{\ge0} \mid \tilde e_i^k\ulambda\neq 0\},\quad
\varphi_i(\ulambda)=\max\{k\in\Z_{\ge0} \mid \tilde f_i^k\ulambda\neq 0\}.
$$

The empty $l$-partition $\boldsymbol\emptyset=(\emptyset,...,\emptyset)$ 
is a highest weight vector of weight $\Lambda$ in $\Fock^{\mathbf v}$. 
We denote by $V_v(\mathbf v)$ the $U_v(\mathfrak g)$-submodule 
generated by $\boldsymbol\emptyset$. Then, 
$V_v(\mathbf v)$ is isomorphic to $V_v(\Lambda)$ as $U_v(\mathfrak g)$-modules. 

\begin{defn}
The crystal $B(\mathbf v)$ is 
the connected subcrystal of $B^{\mathbf v}$ that contains the 
empty $l$-partition $\boldsymbol\emptyset$. An $l$-partition in 
$B(\mathbf v)$ is called 
an \emph{Uglov $l$-partition} of multicharge $\mathbf v$.
\end{defn}

As $B(\mathbf v)$ is the subcrystal which corresponds to $V_v(\mathbf v)$, 
it is isomorphic to the highest weight crystal $B(\Lambda)$. 

\subsection{FLOTW $l$-partitions} 
Define a set $\mathcal V_l$ of multicharges by
$$
\mathcal V_l=
\{\mathbf v=(\mathrm v_0,\ldots,\mathrm v_{l-1}) \mid
\mathrm v_0\leq \mathrm v_1\leq\cdots \leq \mathrm v_{l-1}<\mathrm v_0+e\}.
$$

For each $l$-partition $\ulambda=(\lambda^{(0)},\dots,\lambda^{(l-1)})$, 
let $\lambda_j^{(c)}$, for $j=1,2,\cdots$, be the parts of $\lambda^{(c)}$. 
If $\lambda_j^{(c)}>0$ then 
we denote the residue of the right end node of the $j^{th}$ row 
of $\lambda^{(c)}$ by $\residue(\lambda_j^{(c)})$, which is the residue modulo $e$ of 
$\lambda_j^{(c)}-j+\mathrm v_c$. 

\begin{defn}
Suppose that $\mathbf v\in \mathcal{V}_{l}$. 
A \emph{FLOTW $l$-partition} of multicharge $\mathbf v$ is 
an $l$-partition $\ulambda$ which 
satifies the following two conditions. 

\begin{enumerate}
\item[(i)]
We have the inequalities
$$
\lambda_j^{(c)}\geq\lambda^{(c+1)}_{j+\mathrm v_{c+1}-\mathrm v_c}, 
\;\text{for $0\le c\le l-2$},\;\; \text{and}\;\;
\lambda^{(l-1)}_j\geq\lambda^{(0)}_{j+e+\mathrm v_0-\mathrm v_{l-1}}.
$$
\item[(ii)]
For each $k\in\Z_{>0}$, we have 
$$
\{\residue(\lambda_j^{(c)}) \mid \lambda_j^{(c)}=k\}\neq \Z/e\Z. 
$$
\end{enumerate}
We denote by $\Phi(\mathbf v)_n$ the set of FLOTW $l$-partitions 
of multicharge $\mathbf v$ and rank $n$. Then, we define 
$$
\Phi(\mathbf v)=\bigsqcup_{n\ge0}\Phi(\mathbf v)_n,\quad\text{and}\quad
\Phi=\bigsqcup_{\mathbf v\in\mathcal V_l}\Phi(\mathbf v).
$$
\end{defn}

We have the following result \cite{FLOTW}. 

\begin{prop}
Suppose that $\mathbf v\in\mathcal V_l$. Then, 
$B(\mathbf v)=\Phi(\mathbf v)$.
\end{prop}

\subsection{Kleshchev $l$-partitions} 
If $l=1$ then we have the level $1$ deformed Fock spaces $\Fock^{\mathrm v}$, for 
$\mathrm v\in\Z$. We consider the tensor product
$$
\Fock^{\mathrm v_{l-1}}\otimes\cdots\otimes\Fock^{\mathrm v_0},
$$
for a multicharge $\mathbf v$. 
Note that it depends only on 
$\overline{\mathbf v}=(\overline{\mathrm v}_0,\dots,\overline{\mathrm v}_{l-1})$. 
Then, 
$$
\left(L^{\mathrm v_{l-1}}\otimes\cdots\otimes L^{\mathrm v_0},
B^{\mathrm v_{l-1}}\otimes\cdots\otimes B^{\mathrm v_0}\right)
$$
is a crystal basis of 
$\Fock^{\mathrm v_{l-1}}\otimes\cdots\otimes\Fock^{\mathrm v_0}$. 

\begin{defn}
A \emph{Kleshchev $l$-partition} is 
an $l$-partition $\ulambda$ 
such that the tensor product of the transpose of 
$\lambda^{(i)}$'s in the reversed order 
$$
\hphantom{}^t\lambda^{(l-1)}\otimes\cdots\otimes\hphantom{}^t\lambda^{(0)}
$$
belongs the connected component of 
$B^{\mathrm v_{l-1}}\otimes\cdots\otimes B^{\mathrm v_0}$ that contains 
$\emptyset\otimes\cdots\otimes\emptyset$. 

We denote by $\Phi^K_n$ the set of Kleshchev $l$-partitions of rank $n$. 
Then we define
$$
\Phi^K=\bigsqcup_{n\ge 0}\Phi^K_n.
$$
\end{defn}
We need the transpose of partitions in the definition 
in order to make it compatible with Specht module theory of 
cyclotomic Hecke algebras, which we introduce later. 
Note that if $\ulambda$ is Kleshchev then each component $\lambda^{(j)}$ is 
$e$-restricted. 

$\Phi^K$ inherits the crystal structure from 
$B^{\mathrm v_{l-1}}\otimes\cdots\otimes B^{\mathrm v_0}$, and 
$\Phi^K$ is isomorphic to the highest weight crystal 
$B(\Lambda)$, again. 

\subsection{Crystal isomorphisms} 
As $\Phi(\mathbf v)$ and $\Phi^K$ are isomorphic, 
we have a unique isomorphism of crystals between them, which we denote by
$$
\Gamma :\Phi(\mathbf v)\rightarrow \Phi^K.
\label{gamma}
$$

We may compute this bijection explicitly. In fact, if we fix $n$ and 
choose another multicharge $\mathbf w$ such that 
\begin{itemize}
\item
$\mathrm w_j$ is sufficiently smaller than 
$\mathrm w_{j+1}$, for $0\le j\le l-2$, and 
\item
$\overline{\mathrm v}_j=\overline{\mathrm w}_j$, for $0\le j\le l-1$, 
\end{itemize}
then the bijection 
between $\Phi^K_{\le n}$ and $B(\mathbf w)_{\le n}$ given by 
$$
(\lambda^{(0)},\dots,\lambda^{(l-1)})\mapsto 
(\hphantom{}^t\lambda^{(0)},\dots,\hphantom{}^t\lambda^{(l-1)})
$$
is compatible with the crystal structures on 
$\Phi^K_{\le n}$ and $B(\mathbf w)_{\le n}$. 
Hence, it suffices to compute the crystal isomorphism between 
$B(\mathbf v)$ and $B(\mathbf w)$. 

Let $\widehat{\Sym}_n=e\Z\wr\Sym_n\subseteq\Aut(\Z^l)$ be 
the extended affine symmetric group. Define $\sigma_j\in\Aut(\Z^l)$, 
for $0\le j\le l-2$, by 
$$
\sigma_j(\mathrm v_0,\dots,\mathrm v_{j-1},\mathrm v_j,\dots,\mathrm v_{l-1})
=(\mathrm v_0,\dots,\mathrm v_j,\mathrm v_{j-1},\dots,\mathrm v_{l-1})
$$
and define $\tau\in\Aut(\Z^l)$ by
$\tau(\mathrm v_0,\dots,\mathrm v_{l-1})
=(\mathrm v_1,\dots,\mathrm v_{l-1},\mathrm v_0+e)$. Then, 
$\widehat{\Sym}_n$ is generated by these elements. 

The following theorem was proved by the second and 
the third authors in \cite{JL2}. 
As the multicharges $\mathbf v$ and $\mathbf w$ are in the 
same $\widehat{\Sym}_n$-orbit, it allows us to compute 
the crystal isomorphism between 
$B(\mathbf v)$ and $B(\mathbf w)$ explicitly. 

\begin{thm}
\label{JL-isoms}
\item[(1)]
The crystal isomorphism $B(\mathbf v)\rightarrow B(\tau\mathbf v)$ is 
given by 
$$
(\lambda^{(0)},\dots,\lambda^{(l-1)})\mapsto
(\lambda^{(1)},\dots,\lambda^{(l-1)},\lambda^{(0)}).
$$
\item[(2)]
The crystal isomorphism $B(\mathbf v)\rightarrow B(\sigma_j\mathbf v)$ is 
given by 
$$
(\lambda^{(0)},\dots,\lambda^{(j-1)},\lambda^{(j)},\dots,\lambda^{(l-1)})
\mapsto
(\lambda^{(0)},\dots,\tilde\lambda^{(j)},\tilde\lambda^{(j-1)},\dots,\lambda^{(l-1)}),
$$
where, $\tilde\lambda^{(j-1)}$ and $\tilde\lambda^{(j)}$ are defined by 
$$
\lambda^{(j)}\otimes\lambda^{(j-1)}\mapsto
\tilde\lambda^{(j-1)}\otimes\tilde\lambda^{(j)}
$$
under the following crystal isomorphism, called a combinatorial R-matrix, between
$\mathfrak g(A_\infty)$-crystals. 
$$
B(\Lambda_{\mathrm v_j})\otimes B(\Lambda_{\mathrm v_{j-1}})\rightarrow
B(\Lambda_{\mathrm v_{j-1}})\otimes B(\Lambda_{\mathrm v_j}).
$$
\end{thm}

The combinatorial R-matrix may be computed in a purely combinatorial manner. See 
\cite{JL2} for the details. 

\subsection{Crystal embedding to $B(\infty)$} 
Let $T_\Lambda=\{t_\Lambda\}$ be the crystal defined by 
$\wt(t_\Lambda)=\Lambda$, $\epsilon_i(t_\Lambda)=\varphi_i(t_\Lambda)=-\infty$ 
and $\tilde e_it_\Lambda=\tilde f_it_\Lambda=0$. Then, by the theory 
of crystals, we have the crystal embedding
$B(\Lambda)\hookrightarrow B(\infty)\otimes T_\Lambda$ such that
\begin{itemize}
\item[(i)]
the image of the embedding is given by
$$
\{ b\otimes t_\Lambda\in B(\infty)\otimes T_\Lambda \mid 
\epsilon_i(b^*)\le \Lambda(\alpha_i^\vee)\},
$$
where $b\mapsto b^*$ is the involution on $B(\infty)$ which is 
induced by the anti-automorphism of $U_v^-$ defined by $f_i\mapsto f_i$, 
\item[(ii)]
$b\otimes t_\Lambda$ belongs to the image if and only if  
$G_v(b)v_\Lambda\neq 0$, 
where $v_\Lambda$ is the highest weight vector of $V_v(\Lambda)$. 
\end{itemize}

We identify $B(\infty)$ with the crystal of aperiodic multisegments 
defined in Theorem \ref{LTV-2} and used in 
Theorem \ref{The second branching theorem}. 
As $B(\mathbf v)$ is isomorphic to $B(\Lambda)$, the 
crystal embedding defines a map
$$
B(\mathbf v)\hookrightarrow B(\infty)\otimes T_\Lambda.
$$
We shall describe this map in subsequent subsections. 
By virtue of Theorem \ref{JL-isoms}, we may assume that 
$\mathbf v\in\mathcal V_l$. Write the crystal embedding by 
$\ulambda\mapsto f(\ulambda)\otimes t_\Lambda$, and denote both 
the empty $l$-partition and the empty multisegment by the 
common symbol $\mathbf\emptyset$. 
Then, the crystal embedding  
sends $\mathbf\emptyset$ to $\mathbf{\emptyset}\otimes t_\Lambda$, 
and the tensor product rule shows that for any path 
$$
\mathbf\emptyset\overset{i_1}{\rightarrow}\ulambda_1
\overset{i_2}{\rightarrow}\ulambda_2
\overset{i_3}{\rightarrow}\cdots 
\overset{i_n}{\rightarrow}\ulambda_n
$$
in $B(\mathbf v)$, we have the corresponding path 
$$
\mathbf\emptyset\overset{i_1}{\rightarrow}f(\ulambda_1)
\overset{i_2}{\rightarrow}f(\ulambda_2)
\overset{i_3}{\rightarrow}\cdots
\overset{i_n}{\rightarrow}f(\ulambda_n)
$$
in $B(\infty)$, and vice versa. 
On the other hand, if one can prove this property for some map 
$f: B(\mathbf v)\rightarrow B(\infty)$ then it follows that 
$$
\epsilon_i(\ulambda)=\epsilon_i(f(\ulambda)\otimes t_\Lambda)\;\;
\text{and}\;\;
\wt(\ulambda)=\wt(f(\ulambda)\otimes t_\Lambda),
$$
so that 
we also have $\varphi_i(\ulambda)=\varphi_i(f(\ulambda)\otimes t_\Lambda)$. 
Hence, 
we may conclude that the map $\ulambda\mapsto f(\ulambda)\otimes t_\Lambda$ 
is a crystal embedding in the sense of \cite{Ka} and 
it must coincide with the crystal embedding $B(\mathbf v)\hookrightarrow 
B(\infty)\otimes T_\Lambda$. 

\subsection{Row lengths and the order $\prec_{\mathbf v}$} 
We prove two lemmas which relate the length of rows of an 
$l$-partition and the order $\prec_{\mathbf v}$. 

\begin{lem}
\label{order1} 
Let $\mathbf v\in \mathcal V_l$ and 
$\ulambda=(\lambda^{(0)},\dots,\lambda^{(l-1)})
\in \Phi(\mathbf v)$. Suppose that $\gamma_1=(a_1,b_1,c_1)$ and 
$\gamma_2=(a_2,b_2,c_2)$ are $i$-nodes of $\ulambda$ such that 
each node is either addable or removale $i$-node. 
Then, $\lambda_{a_1}^{(c_1)}<\lambda_{a_2}^{(c_2)}$ implies 
$\gamma_1\prec_{\mathbf v}\gamma_2$.
\end{lem}
\begin{proof}
We show that $\gamma_2\preceq_{\mathbf v}\gamma_1$ implies 
$\lambda_{a_1}^{(c_1)}\geq \lambda_{a_2}^{(c_2)}$. As an intermediate step, 
we first claim that 
$\gamma_2\preceq_{\mathbf v}\gamma_1$ implies 
$\lambda_{a_1}^{(c_1)}\geq \lambda_{b_1-b_2+a_2}^{(c_2)}$. Note that 
we have $c(\gamma_1)\geq c(\gamma_2)$ 
by $\gamma_2\preceq_{\mathbf v}\gamma_1$. Hence, we have
$$
a_1\leq b_1-b_2+a_2+\mathrm v_{c_1}-\mathrm v_{c_2},
$$
which implies 
$\lambda_{a_1}^{(c_1)}\geq 
\lambda_{b_1-b_2+a_2+\mathrm v_{c_1}-\mathrm v_{c_2}}^{(c_1)}$. 

Suppose that $c_1\leq c_2$. As $\ulambda$ is a FLOTW $l$-partition, we have
$$
\lambda_{b_1-b_2+a_2+\mathrm v_{c_1}-\mathrm v_{c_2}}^{(c_1)}
\geq
\lambda_{b_1-b_2+a_2+\mathrm v_{c_1+1}-\mathrm v_{c_2}}^{(c_1+1)}
\geq\cdots\geq
\lambda_{b_1-b_2+a_2}^{(c_2)}.
$$
Hence $\lambda_{a_1}^{(c_1)}\geq \lambda_{b_1-b_2+a_2}^{(c_2)}$ follows. 

Suppose that $c_1>c_2$. Then, $c(\gamma_1)>c(\gamma_2)$ and we must
have 
$$
b_1-a_1+\mathrm v_{c_1}\geq b_2-a_2+\mathrm v_{c_2}+e,
$$
because $\gamma _{1}$ and $\gamma _{2}$ have the same residue $i$. 
Hence, we have
$$
a_1\leq b_1-b_2+a_2+\mathrm v_{c_1}-\mathrm v_{c_2}-e,
$$
which implies 
$\lambda_{a_1}^{(c_1)}\geq 
\lambda_{b_1-b_2+a_2+\mathrm v_{c_1}-\mathrm v_{c_2}-e}^{(c_1)}$. 
Then, by the same reasoning as above, we have
\begin{multline*}
\lambda_{b_1-b_2+a_2+\mathrm v_{c_1}-\mathrm v_{c_2}-e}^{(c_1)}
\geq
\lambda_{b_1-b_2+a_2+\mathrm v_{c_1+1}-\mathrm v_{c_2}-e}^{(c_1+1)}
\geq\cdots\\
\geq
\lambda_{b_1-b_2+a_2+\mathrm v_{l-1}-\mathrm v_{c_2}-e}^{(l-1)}
\geq
\lambda_{b_1-b_2+a_2+\mathrm v_0-\mathrm v_{c_2}}^{(0)}
\geq\cdots\geq
\lambda_{b_1-b_2+a_2}^{(c_2)}.
\end{multline*}
Hence $\lambda_{a_1}^{(c_1)}\geq \lambda_{b_1-b_2+a_2}^{(c_2)}$ follows again. 

If $b_1\leq b_2$ then $b_1-b_2+a_2\leq a_2$ implies the desired inequality  
$\lambda_{a_1}^{(c_1)}\geq \lambda_{a_2}^{(c_2)}$. Suppose that $b_1>b_2$. 
As $\gamma_1$ is either addable or removable $i$-node, 
we have either  
$b_1=\lambda_{a_1}^{(c_1)}+1$ or $b_1=\lambda_{a_1}^{(c_1)}$. Similarly, 
we have either $b_2=\lambda_{a_2}^{(c_2)}+1$ or $b_2=\lambda_{a_2}^{(c_2)}$. 
Hence, we have $\lambda_{a_1}^{(c_1)}\geq b_1-1\geq b_2\geq\lambda_{a_2}^{(c_2)}$. 
\end{proof}

\begin{lem}
\label{order2}
Let $\ulambda$ be a FLOTW $l$-partition, and let 
$\gamma_A=(a',b+1,c')$ and $\gamma_R=(a,b,c)$ be addable and 
removable $i$-nodes of $\ulambda$ respectively. 
Then we have $\gamma_R\prec_{\mathbf v}\gamma_A$.
\end{lem}
\begin{proof}
Suppose to the contrary that $\gamma_A\prec_{\mathbf v}\gamma_R$. Then 
we have either
\begin{enumerate}
\item[(i)]
$c(\gamma_A)<c(\gamma_R)$, or
\item[(ii)]
$c(\gamma_A)=c(\gamma_R)$ and $c'>c$.
\end{enumerate}

In case (i), $b-a+\mathrm v_c\geq b+1-a'+\mathrm v_{c'}+e$ 
so that $a+\mathrm v_{c'}-\mathrm v_c+e\leq a'-1$. 
As $\gamma_A$ is an addable node, we also have 
$\lambda_{a'-1}^{(c')}>\lambda_{a'}^{(c')}$. Then,  
$a+\mathrm v_{c'}-\mathrm v_c+e<a'$ implies that 
$$
\lambda_{a+\mathrm v_{c'}-\mathrm v_c}^{(c')}
\geq
\lambda_{a+\mathrm v_{c'}-\mathrm v_c+e}^{(c')}>\lambda_{a'}^{(c')}.
$$
Now, using the assumption that $\ulambda$ is a FLOTW $l$-partition, 
we have
$$
\begin{cases}
\lambda_a^{(c)}\geq \lambda_{a+\mathrm v_{c'}-\mathrm v_c}^{(c')}
>\lambda_{a'}^{(c')} & \text{ if }c\leq c', \\ 
\lambda_a^{(c)}\geq \lambda_{a+\mathrm v_{c'}-\mathrm v_c+e}^{(c')}
>\lambda_{a'}^{(c')} & \text{ if }c>c'.
\end{cases}
$$
However, $\lambda_a^{(c)}=b$ since $\gamma_R$ is a removable node, and 
$\lambda_{a'}^{(c')}=b$ since $\gamma_A$ is an addable node. 
Thus, we have reached a contradiction. 

In case (ii), $b-a+\mathrm v_c=b+1-a'+\mathrm v_{c'}$ implies 
$a+\mathrm v_{c'}-\mathrm v_c+1=a'$. As $\gamma_A$ is an addable 
node, $\lambda^{(c')}_{a+\mathrm v_{c'}-\mathrm v_c}>\lambda^{(c')}_{a'}$. 
Thus, $c'>c$ implies that
$$
\lambda_a^{(c)}\geq \lambda_{a+\mathrm v_{c'}-\mathrm v_c}^{(c')}
>\lambda_{a'}^{(c')}. 
$$
However, we have $\lambda_a^{(c)}=b$ and $\lambda_{a'}^{(c')}=b$ as before, so that 
we have reached a contradiction again. 
\end{proof}

\subsection{The map $f_{\mathbf v}$} 
For each FLOTW $l$-partition $\ulambda\in \Phi(\mathbf v)$, 
we associate a multisegment which 
is a collection of segments 
$$
[1-i+\mathrm v_c;\lambda^{(c)}_i), 
$$
where $\lambda^{(c)}_i$ are parts of $\lambda^{(c)}$, for 
$c=0,\dots,l-1$. This defines a well-defined map
$f_{\mathbf v}: \Phi(\mathbf v)\rightarrow B(\infty)$. 

\begin{exmp}
Let $e=4$, and let $\ulambda=((2,1),(1))\in\Phi((0,1))$. Then 
$$
f_{(0,1)}(\ulambda)=\{[0,1], [3], [1]\}. 
$$
Next let $\ulambda=((2),(1),(1))\in\Phi((0,1,3))$. Then we have the same
result
$$
f_{(0,1,3)}(\ulambda)=\{[0,1], [1], [3]\}.
$$
\end{exmp}

Then we may prove the following. Note that the fact itself was observed 
by several people including the first author years ago, 
but the authors do not know any reference which proves this. 

\begin{thm}
\label{Th_fv}
Suppose that $\mathbf v\in \mathcal V_l$. Then, 
the crystal embedding $\Phi(\mathbf v)\hookrightarrow 
B(\infty)\otimes T_\Lambda$ is given by 
$\ulambda\mapsto f_{\mathbf v}(\ulambda)\otimes t_\Lambda$. 
\end{thm}
\begin{proof}
As was explained in the previous subsection, 
it suffices to show that there is an arrow 
$$
\ulambda\overset{i}{\rightarrow} \umu
$$ 
in $B(\mathbf v)$ if and only if there is an arrow
$$
f_{\mathbf v}(\ulambda)\overset{i}{\rightarrow} f_{\mathbf v}(\umu)
$$ 
in $B(\infty)$. 

We read the addable and removable $i$-nodes of $\umu$ 
in increasing order with respect to the total order $\prec_{\mathbf v}$.  
Let $\gamma_1\dots \gamma_m$ be the resulting word of the nodes. 
On the other hand, we read the same set of addable and removable $i$-nodes of 
$\umu$ in increasing order with respect to the length of 
the corresponding segments in $f_{\mathbf v}(\umu)$. 
If the length are the same, we declare 
that removable $i$-nodes precede addable $i$-nodes. We denote 
the resulting word $\gamma_{\sigma(1)}\dots\gamma_{\sigma(m)}$, 
for $\sigma\in\Sym_m$. 

Write $\gamma_j=(a_j,b_j,c_j)$, for $1\le j\le m$. Then, 
Lemma \ref{order1} implies that if 
$\lambda_{a_j}^{(c_j)}\neq\lambda_{a_k}^{(c_k)}$ then 
$j<k$ implies $\sigma^{-1}(j)<\sigma^{-1}(k)$. 
On the other hand, Lemma \ref{order2} implies that 
if $\lambda_{a_j}^{(c_j)}=\lambda_{a_k}^{(c_k)}$ then 
$j<k$ implies $\sigma^{-1}(j)<\sigma^{-1}(k)$. We conclude 
that $\sigma$ is the identity. 

We define $S'_{k,i}$ to be the number of addable $i$-nodes minus the number of 
removable $i$-nodes in $\{\gamma_k,\gamma_{k+1},\dots,\gamma_m\}$. 

Suppose that $\tilde e_i\umu=\ulambda$ and 
let $\gamma=(a,b,c)$ be the good $i$-node of $\umu$. Then  
$\min_{k>0} S'_{k,i}$ is attained at $\gamma$. Define $k_r$, for $r>0$, by
$$
k_r=\min\{ j \mid \lambda_{a_j}^{(c_j)}\geq r\}. 
$$
It is clear that $\min_{k>0} S'_{k,i}$ is attained only at 
removable nodes of the form $\gamma_{k_r}$, for some $r$. 
Now observe that addable and removable $i$-nodes of the multisegment 
$f_{\mathbf v}(\umu)$ which do not belong to $\{\gamma_1,\dots,\gamma_m\}$ 
come from pairs of consecutive rows of the same length in $\umu$. 
Let $m_{(k;i]}$ be the multiplicity of $(k;i]$ in $f_{\mathbf v}(\umu)$. 
Then, by the above observation, we have 
$$
S_{r,i}=\sum_{k\geq r}(m_{(k,i-1]}-m_{(k;i]})=S'_{k_r,i},
$$
and $\min_{r>0} S_{r,i}$ is attained at $r=b$. Instead of proving that 
$b$ is the unique $r$ that attains the minimum, we shall show that 
$\tilde f_if_{\mathbf v}(\ulambda)=f_{\mathbf v}(\umu)$. 
As $\gamma$ is the good removable $i$-node of $\umu$, the following is clear. 
\begin{quote}
If $r<b$ then, among the nodes $\gamma_j$, for $k_r\le j<k_b$, 
the number of addable nodes is always greater than or equal 
to the number of removable nodes. 
\end{quote}
This implies that, if we change the status of $\gamma$ from a removable node 
to an addable node, then 
$S_{r,i}>S_{b,i}$ if $r\leq b-1$, 
for the new values $S_{r,i}$ and $S_{b,i}$ 
computed after we change the status of $\gamma$. 
If we consider normal $i$-nodes which appear to the right of $\gamma$, 
it is also clear that $S_{r,i}\geq S_{b,i}$ if $r\geq b+1$, for the new values 
$S_{r,i}$. Thus, we obtain 
$\tilde f_if_{\mathbf v}(\ulambda)=f_{\mathbf v}(\umu)$. 

Next suppose that $\tilde f_if_{\mathbf v}(\ulambda)=f_{\mathbf v}(\umu)$. 
We consider $S_{r,i}$ and suppose that $\min_{r>0} S_{r,i}$ is attained 
at $\ell_0<\ell_1<\cdots$. The minimum value is attained at a removable $i$-node
which is the leftmost node among the nodes of the segments of the same length. 
Then, the minimality implies that the right neighbor of the removable node is 
addable. We denote this node by $\gamma$. 
We show that $\gamma$ is the good addable $i$-node of $\ulambda$. 

Suppose that $\gamma$ is cancelled in the RA-deletion procedure. 
If the removable $i$-node which cancels $R$ is not of the form $\gamma_{k_r}$, 
it contradicts the minimality of $S_{\ell_0,i}$. Thus, the removable node is 
$\gamma_{k_b}$, for some $b<\ell_0$. Then, $S_{b,i}=S_{\ell_0,i}$ implies 
$\ell_0\leq b$, which contradicts $b<\ell_0$. Hence, we have proved that 
$\gamma$ is a normal addable $i$-node. If there was another normal addable $i$-node 
to the right of $\gamma$, it would contradict the minimality of $S_{\ell_0,i}$, 
so that $\gamma$ is the good addable $i$-node of $\ulambda$. Thus, we obtain 
$\tilde f_i\ulambda=\umu$. 
\end{proof}

Define $B^{\mathrm ap}(\Lambda)=\{ \psi\in B(\infty) \mid 
\epsilon_i(\psi^*)\le \Lambda(\alpha_i^\vee)\}$. 
As we have proved that $\ulambda\mapsto f_{\mathbf v}(\ulambda)\otimes t_{\Lambda}$ 
is the crystal embedding $B(\Lambda)\hookrightarrow B(\infty)\otimes T_\Lambda$ 
in the language of FLOTW and multisegment realizations, 
we have the following corollary. The basis in Corollary \ref{Cor_ap}(2) is the 
\emph{canonical basis} of $V_v(\Lambda)$. The statement is for the crystal 
structure we have chosen, 
but it is easy to state it for the other, since the Kashiwara involution 
on the set of aperiodic multisegments is explicitly described in \cite{JL1}. 

\begin{cor}
\label{Cor_ap}
\item[(1)]
$f_{\mathbf v}(B(\mathbf v))=B^{\mathrm{ap}}(\Lambda)$.
\item[(2)]
$\{ G_v(\psi)v_\Lambda \mid \psi\in B^{\mathrm ap}(\mathrm v)\}$ is 
a basis of $V_v(\Lambda)$. 
\end{cor}

\section{Fock space theory for cyclotomic Hecke algebras}

In this section, we give the combinatorial proof of the modular branching rule. 
The proof depends on Lemma \ref{comparison-2}, which 
says that isomorphisms of crystals give the correspondence of labels of a simple 
$\mathcal H_n^\Lambda$-module, which is labelled by various realizations 
of the crystal $B(\Lambda)$. 
Hence, the explicit description of the isomorphisms 
in the previous section gives us the module correspondence. 

\subsection{Cyclotomic Hecke algebras} 
Let $\mathbf v$ be a multicharge as before. 
The \emph{cyclotomic Hecke algebra} $\mathcal H_n^{\mathbf v}(q)$ is the 
quotient algebra $\Hecke_n/I_{\mathbf v}$ of the affine Hecke algebra $\Hecke_n$, where 
$I_{\mathbf v}$ is the ideal of $\Hecke_n$ 
generated by the polynomial $\prod_{i=0}^{l-1}(X_1-q^{\mathrm v_i})$. 
If we specialize $q=\zeta$, the algebra depends only on $\Lambda$, and 
we denote the algebra by $\mathcal H_n^\Lambda$. This is the main object 
of the study in the remaining part of the paper. As $\mathcal H_n^\Lambda$ 
is a quotient algebra of the affine Hecke algebra $\Hecke_n$, the set of 
simple $\mathcal H_n^\Lambda$-modules is a subset of simple $\Hecke_n$-modules. 
In fact, by Fock space theory for cyclotomic Hecke algebras 
we will explain in the next subsection, 
we know that it is the set $\{ D_\psi \mid \psi \in B^{\mathrm ap}(\Lambda)\}$.

\begin{defn}
We denote by $\mathcal H_n^\Lambda\modA$ the category of
finite-dimensional $\mathcal H_n^\Lambda$-modules. 
\end{defn}

Note that $\mathcal H_n^{\mathbf v}(q)$ is a cellular algebra in the sense of 
Graham and Lehrer: it has the Specht module theory developped by Dipper, James 
and Mathas. Then, the first author showed that  
simple $\mathcal H_n^\Lambda$-modules are labelled by Kleshchev $l$-partitions. 
We refer to \cite[Ch. 12]{Ari3} for details.

For $\ulambda\in \Phi^K$, we denote by 
$D^\ulambda$ the simple $\mathcal H_n^\Lambda$-module labelled 
by $\ulambda$. For $\ulambda\in \Phi(\mathbf v)_n$, we define 
$\widetilde D^\ulambda$ by
$$
\widetilde D^\ulambda=D^{\Gamma(\ulambda)}.  \label{D-tilde}
$$
We will explain in the next subsection that 
this labelling coincides with the Geck-Rouquier-Jacon parametrization of 
simple $\mathcal H_n^\Lambda$-modules in terms of the canonical basic set. 

Before giving the second proof, we complete the first proof. Namely, 
we prove Theorem \ref{comparison-1} below, 
which compares the geometrically defined simple 
$\mathcal H_n^\Lambda$-modules and the combinatorially defined simple 
$\mathcal H_n^\Lambda$-modules by using Theorem \ref{The second branching theorem}. 

\begin{thm}
\label{comparison-1}
Let $\ulambda$ be an $l$-partition. Then,   
$\widetilde D^\ulambda\simeq D_{f_{\mathbf v}(\ulambda)}$ as $\Hecke_n$-modules.
\end{thm}
\begin{proof}
We have $i\text{-}\Res(D_\psi)\simeq D_{\tilde e_i\psi}$ by 
Theorem \ref{The second branching theorem}. On the other hand, 
we have $i\text{-}\Res(D^\ulambda)\simeq D^{\tilde e_i\ulambda}$, 
for $\ulambda\in\Phi^K$, in 
\cite[Theorem 6.1]{Ari4}. Note that if 
$i\text{-}\Res(D^\ulambda)\simeq i\text{-}\Res(D_\psi)\neq 0$ then 
$D^\ulambda\simeq D_\psi$. This property of crystals is a 
consequence of the Frobenius reciprocity. Hence, we may 
prove the claim by induction on $n$. 
\end{proof}

\subsection{Standard modules}
We say a few words on the standard modules of the affine Hecke algebra. 
Let $X\in \mathcal O_\psi$ and consider 
$$
(\Fl_n^a)_X=\{F\in \Fl_n^a \mid XF_i\subseteq F_{i-1}\}.
$$ 
Then, $H_*((\Fl_n^a)_X,\C)$ is an $H_*^{BM}(Z_n^a,\C)$-module 
by the convolution action, and it is called the \emph{standard module}. 
We denote it by $M_\psi$. 
Suppose that $X$ is a principal nilpotent element so that 
$\psi=[i;l)$ for some $i\in\Z/e\Z$ and $l\in\Z_{>0}$. Then, 
$(\Fl_n^a)_X$ is a point, which is the flag 
$$
0\subseteq \Ker(X)\subseteq \Ker(X^2)\cdots\subseteq \Ker(X^n)=V
$$
of flag type $(i+l-1,\dots,i+1,i)$, and the proof of 
Lemma \ref{eigenvalue analysis} shows that, 
if we follow the identification 
$\Hecke_n\simeq K^{G_n\times\C^\times}(Z_n)$ in \cite{L3}, 
then $M_\psi$ is the one dimensional 
$\Hecke_n$-module given by $T_i\mapsto -1$ and 
$$
X_1\mapsto \zeta^{i+l-1},\dots, X_{n-1}\mapsto \zeta^{i+1}, X_n\mapsto \zeta^i. 
$$
Thus, $M_\psi$ for general $\psi$ coincides with the induced up module 
of the tensor product of such one dimensional modules over the 
affine Hecke algebras associated with segments in $\psi$, in the Grothendieck 
group of the module category of the affine Hecke algebra. 

Now, we switch to the other identification used in 
Theorem \ref{The second branching theorem}. Define the \emph{standard module} $N_\psi$ by 
$$
N_\psi=\hphantom{}^\sigma M_{\rho(\psi)}.
$$
Then $N_\psi$ is given by $T_i\mapsto \zeta$ and 
$$
X_1\mapsto \zeta^{i-l+1},\dots, X_{n-1}\mapsto \zeta^{i-1}, X_n\mapsto \zeta^i. 
$$
when $\psi=(l;i]$. This is the standard module in \cite{Ari3}. 
Then, a key observation used in \cite{Ari3} was the equality
$$
G_{v=1}(\psi)=\sum_{\psi'} [N_{\psi'}: D_\psi]u_{\psi'}
$$
in the Hall algebra in subsection 2.3 evaluated at $v=1$. 
\footnote{In fact, the choice of the identification played no role in \cite{Ari3} 
because it sufficed for us to prove the statement that 
the canonical basis evaluated at $v=1$ 
coincides with the dual basis of simples in the Fock space, and we did not need 
compare individual simple modules.} 
Now we are able to give an example of Theorem \ref{comparison-1}. 

\begin{exmp}
Let $e=3$. Then, we have
$$
G_{v=1}(\{(2;2]\})=u_{\{(2;2]\}}+u_{\{(1;1], (1;2]\}},\;\;
G_{v=1}(\{(1;1], (1;2]\})=u_{\{(1;1], (1;2]\}}.
$$
Note that $N_{\{(1;1]\}}$ and $N_{\{(1;2]\}}$ 
are one dimensional $\Hecke_1$-modules defined by $X_1\mapsto \zeta$ and 
$X_1\mapsto \zeta^2$, respectively. Then, $N_{\{(2;2]\}}=D_{\{(2;2]\}}$ is 
the simple module defined by
$$
X_1\mapsto \zeta,\;\; X_2\mapsto \zeta^2, \;\; T_1\mapsto \zeta,
$$ 
and $N_{\{(1;1], (1;2]\}}$ is the module induced from 
$N_{\{(1;1]\}}\otimes N_{\{(1;2]\}}$. Thus, we deduce that 
$D_{\{(1;1], (1;2]\}}$ is the simple module defined by
$$
X_1\mapsto \zeta^2,\;\; X_2\mapsto \zeta,\;\; T_1\mapsto -1. 
$$

\begin{itemize}
\item[(Ex.1)]
Suppose that $l=1$. Then we have
$$
D^{(2)}\simeq D_{\{(2;2]\}},\;\text{for $\mathbf v=1$,} \quad
D^{(1^{2})}\simeq D_{\{(1;1],(1;2]\}}, \;\text{for $\mathbf v=2$,}
$$
with $(2), (1^2)\in \Phi^K_2$. This follows from the explicit construction of 
Specht modules. Since $(2)=\tilde f_2\tilde f_1\emptyset$ and 
$(1^2)=\tilde f_1\tilde f_2\emptyset$ in $\Phi^K_2$, we have 
$\Gamma((2))=(2)$ and $\Gamma((1^2))=(1^2)$, so that
$$
\widetilde D^{(2)}\simeq D_{\{(2;2]\}}\;\; \text{and}\;\; 
\widetilde D^{(1^2)}\simeq D_{\{(1;1], (1;2]\}}. 
$$
\item[(Ex.2)]
Suppose that $l=2$ and $\mathbf v=(1,2)$. Then, 
$\tilde f_2\tilde f_1\mathbf{\emptyset}=((2),\emptyset)$ and 
$\tilde f_1\tilde f_2\mathbf{\emptyset}=((1),(1))$ in $\Phi(\mathbf v)$, 
so that 
$$
\widetilde D^{((2),\emptyset)}\simeq D_{\{(2;2]\}}\;\; \text{and}\;\; 
\widetilde D^{((1),(1))}\simeq D_{\{(1;1], (1;2]\}}. 
$$
\end{itemize}
\end{exmp} 

\subsection{Fock space theory} 

In this subsection, we explain the Fock space theory for cyclotomic Hecke 
algebras. In the following, $G_v(b)$, $U_v^-$, etc. at $v=1$ are denoted by 
$G(b)$, $U^-$, etc. 

Let $\mathcal C_n$ be the full subcategory of $\Hecke_n\modA$ consisting of 
finite dimensional $\Hecke_n$-modules on which $X_1,\dots, X_n$ have eigenvalues in 
$\{1,\zeta,\dots,\zeta^{e-1}\}$. 

\begin{defn}
Let 
$$
U_n=\Hom_\C(K_0(\mathcal C_n),\C)\;\;\text{and}\;\;
V_n=\Hom_\C(K_0(\mathcal H_n^\Lambda\modA),\C)
$$
be the dual spaces of the Grothendieck groups of 
$\mathcal C_n$ and $\mathcal H_n^\Lambda\modA$, and define
$$
U=\bigoplus_{n\geq 0}U_n\;\;\text{and}\;\;
V=\bigoplus_{n\geq 0}V_n.
$$
\end{defn}

Hereafter, we identify $V_n$ with the split Grothendieck group of the additive 
subcategory of $\mathcal H_n^\Lambda\modA$ consisting of 
projective $\mathcal H_n^\Lambda$-modules. 

$U_n$ has the dual basis 
$$
\{[D_\psi]^* \mid \text{$\psi$ is an aperiodic multisegment of rank $n$.}\}
$$ 
which is dual to the basis consisting of simple $\mathcal H_n^\Lambda$-modules. 

Let $\pi: U\rightarrow V$ be the natural map and define 
$$
p: U^-\rightarrow V(\Lambda)\subseteq \Fock
$$
by $F\mapsto Fv_\Lambda$, for $F\in U^-$. 

The theorem below states the most basic result in the Fock space theory. 
See \cite{Ari2} or \cite[Theorem 14.49]{Ari3}.

\begin{thm}
\label{TH_commute}
\item[(1)]
$U$ has structure of a $U^-$-module and $V$ has structure of a $\mathfrak g$-module. 
\item[(2)]
$U$ is isomorphic to the regular representation of $U^-$ such that  
$$
[D_\psi]^*\mapsto G(\psi).
$$
\item[(3)]
$V$ is isomorphic to $V(\Lambda)$ and the basis 
$$
\bigsqcup_{n\geq 0}\{[P] \mid 
\text{$P$ is an indecomposable $\mathcal H_n^\Lambda$-module.}\}
$$
of $V$ corresponds to the canonical basis of $V(\Lambda)$ under the isomorphism. 
\item[(4)]
The following diagram commutes: 
\begin{equation*}
\begin{array}{ccc}
U & \simeq & U^- \\ 
\pi\downarrow\quad &  & \downarrow p \\ 
V & \simeq & V(\Lambda)
\end{array}
\end{equation*}
\end{thm}

\subsection{The combinatorial proof} 
First we make it clear what we mean by 
\lq\lq simple $\mathcal H_n^\Lambda$-modules 
are labelled by Uglov $l$-partitions\rq\rq. 

\begin{defn}
We say that 
\emph{simple $\mathcal H_n^\Lambda$-modules are labelled by $B(\mathbf v)$}, 
if the projective cover of a simple $\mathcal H_n^\Lambda$-module 
is equal to $G(\ulambda)\in \Fock^{\mathbf v}$ in Theorem \ref{TH_commute}(3), 
for $\ulambda\in B(\mathbf v)$, then the label of the simple module is $\ulambda$. 
\end{defn}

It is proved by the first author that Specht module theory is an example of 
the statement that simple 
$\mathcal H_n^\Lambda$-modules are labelled by $B(\mathbf v)$. 
Another example is provided by the second author. Recall that Geck and Rouquier 
invented different theory to label simple modules by using Lusztig's $a$-values. 
The labelling set is called the \emph{canonical basic set}. 
When we work with Hecke algebras of type B, it provides us with a set of 
bipartitions. The second author has generalized the theory to cyclotomic 
Hecke algebras and his result says that simple $\mathcal H_n^\Lambda$-modules 
are labelled by $\Phi(\mathbf v)$, for $\mathbf v\in\mathcal V_l$. 

If one uses Theorem \ref{TH_commute}, it is quite easy to identify 
simple $\mathcal H_n^\Lambda$-modules in various labellings.  

\begin{lem}
\label{comparison-2}
\item[(1)]
Suppose that simple $\mathcal H_n^\Lambda$-modules are labelled by 
$B(\mathbf v)$. Let 
$$
f_{\mathbf v, \infty}: B(\mathbf v)\simeq B^{\mathrm ap}(\Lambda)\subseteq B(\infty)
$$
be the unique crystal isomorphism. Then, $D^\ulambda\simeq D_{f_{\mathbf v}(\ulambda)}$ 
as $\Hecke_n$-modules. 
\item[(2)]
For two labelling $B(\mathbf v)$ and $B(\mathbf w)$ of 
simple $\mathcal H_n^\Lambda$-modules, we denote the set of simple modules by 
$$
\{D_{\mathbf v}^\ulambda \mid \ulambda\in B(\mathbf v)\}\;\;\text{and}\;\;
\{D_{\mathbf w}^\ulambda \mid \ulambda\in B(\mathbf v)\},
$$
respectively. 
Let $f_{\mathbf v, \mathbf w}: B(\mathbf v)\simeq B(\mathbf w)$ be the 
unique crystal isomorphism. Then, 
$D_{\mathbf v}^\ulambda\simeq D_{\mathbf w}^{f_{\mathbf v, \mathbf w}(\ulambda)}$ 
as $\Hecke_n$-modules. 
\end{lem}
\begin{proof}
(1) Suppose that $f_{\mathbf v, \infty}(\ulambda)=\psi$. Then, 
we have $G_v(\psi)\mathbf{\emptyset}=G_v(\ulambda)$. Specializing at 
$v=1$, we obtain $G(\psi)=P^\ulambda$. 
Then, using the commutativity of the diagram in 
Theorem \ref{TH_commute}(4), we conclude that 
$\pi([D_\psi]^*)=[D^\ulambda]^*$, which is identified with $P^\ulambda$. 
Hence, $D_\psi\simeq D^\ulambda$ as $\Hecke_n$-modules. 
 
(2) First we apply (1) to two crystal isomorphisms 
$B(\mathbf v)\simeq B^{\mathrm ap}(\Lambda)$ and 
$B^{\mathrm ap}(\Lambda)\simeq B(\mathbf w)$. Then use the fact that 
$f_{\mathbf v, \mathbf w}=f_{\mathbf w, \infty}^{-1}\circ f_{\mathbf v, \infty}$.
\end{proof}

As we have established Lemma \ref{comparison-2}, we 
can derive the modular branching rule for the affine Hecke algebra from this. 

\begin{thm}
For each aperiodic multisegment $\psi$, we have 
$$
\Soc(i\text{-}\Res_{\Hecke_{n-1}}^{\Hecke_n}(D_\psi))
\simeq D_{\tilde e_i\psi}. 
$$
\end{thm}
\begin{proof}
Choose $\Lambda$ sufficiently large so that 
$f_{\mathbf v}(B(\mathbf v))=B^{\mathrm ap}(\Lambda)$ 
may contain any path 
$$
\mathbf{\emptyset}\overset{i_1}{\rightarrow} \psi_1
\overset{i_2}{\rightarrow} \psi_2
\overset{i_3}{\rightarrow} \cdots
\overset{i_n}{\rightarrow} \psi_n=\psi
$$
in $B(\infty)$ from $\mathbf{\emptyset}$ to $\psi$. 
Let $i\in\Z/e\Z$ be such that $\tilde e_i\psi\neq 0$ and 
let $\ulambda\in B(\mathbf v)$ be such that 
$f_{\mathbf v}(\ulambda)=\psi$. 
Then $\tilde e_i\ulambda\neq 0$ and 
$f_{\mathbf v}(\tilde e_i\ulambda)=\tilde e_i\psi$. 
Then, the previous Lemma yields the isomorphisms 
$$
D_\psi\simeq \widetilde D^\ulambda\;\;\text{and}\;\;
D_{\tilde e_i\psi}\simeq \widetilde D^{\tilde e_i\ulambda}.
$$
Thus, 
$$
\Soc(i\text{-}\Res_{\Hecke_{n-1}}^{\Hecke_n}(D_\psi))
\simeq 
\Soc(i\text{-}\Res_{\Hecke_{n-1}}^{\Hecke_n}(\widetilde D^\ulambda))
\simeq 
\widetilde D^{\tilde e_i\ulambda}
\simeq D_{\tilde e_i\psi},
$$
where the middle isomorphism is the modular branching rule in the labelling 
by Kleshchev $l$-partitions \cite[Theorem 6.1]{Ari4}. 
We have proved the theorem.
\end{proof}

\end{document}